\documentclass[11pt]{article}

\usepackage[utf8]{inputenc}
\usepackage[T1]{fontenc}
\usepackage{lmodern}
\usepackage{microtype}
\usepackage[section]{placeins}
\usepackage{tikz}
\usetikzlibrary{arrows.meta,positioning,shadows.blur,calc}
\usepackage{cite}
\usepackage{graphicx}      
\usepackage{caption}       
\usepackage{subcaption}    
\usepackage{geometry}
\usepackage{amsmath,amssymb}
\usepackage[numbers]{natbib}
\usepackage{doi}
\usepackage{enumitem}
\setlist{nosep}

\usepackage{graphicx}
\usepackage{xcolor}
\usepackage{hyperref}
\hypersetup{
  colorlinks=true,
  linkcolor=blue,
  urlcolor=blue,
  citecolor=blue
}

\usepackage{authblk}

\title{Interpretable AI-Assisted Early Reliability Prediction for a Two-Parameter Parallel Root-Finding Scheme}
\author[4]{Bruno Carpentieri\thanks{Correspondence: bruno.carpentieri@unibz.it}}
\author[3]{Andrei Velichko}
\author[1,2]{Mudassir Shams}
\author[4]{Paola Lecca}
\affil[1]{Department of Mathematics, Faculty of Arts and Science, Balikesir University, 10145 Bal\i kesir, Turkey; mudassir.shams@balikesir.edu.tr}
\affil[2]{Department of Mathematics and Statistics, Riphah International University, I-14, Islamabad 44000, Pakistan}
\affil[3]{Institute of Physics and Technology, Petrozavodsk State University, 185910 Petrozavodsk, Russia; velichko@petrsu.ru}
\affil[4]{Faculty of Engineering, Free University of Bozen--Bolzano, 39100 Bolzano, Italy}
\date{}

\begin{document}
\maketitle

\begin{abstract}
We propose an interpretable AI-assisted reliability diagnostic framework for parameterized root-finding schemes based on kNN–LLE proxy stability profiling and multi-horizon early prediction. The approach augments a numerical solver with a lightweight predictive layer that estimates solver reliability from short prefixes of iteration dynamics, enabling early identification of stable and unstable parameter regimes. For each configuration in the parameter space, raw and smoothed proxy profiles of a largest Lyapunov exponent (LLE) estimator are constructed, from which contractivity-based reliability scores summarizing finite-time convergence are derived. Machine learning models predict the reliability score from early segments of the proxy profile, allowing the framework to determine when solver dynamics become diagnostically informative. Experiments on a two-parameter parallel root-finding scheme show reliable prediction after only a few iterations: the best models achieve $R^2 \approx 0.48$ at horizon $T=1$, improve to $R^2 \approx 0.67$ by $T=3$, and exceed $R^2 \approx 0.89$ before the characteristic minimum-location scale of the stability profile. Prediction accuracy increases to $R^2 \approx 0.96$ at larger horizons, with mean absolute errors around $0.03$, while inference costs remain negligible (microseconds per sample). The framework provides interpretable stability indicators and supports early decisions during solver execution, such as continuing, restarting, or adjusting parameters.
\end{abstract}


\section{Introduction}
\label{sec:introduction}

Robust and reliable root-finding remains a central building block in scientific computing, underpinning nonlinear optimization, inverse problems, parameter estimation, and the solution of discretized PDE models \cite{ortega1970iterative,kelley2003solving,dennis1996numerical}. In many applications, the practical performance of an iterative solver depends not only on the underlying function $f$ and the initialization, but also on algorithmic hyperparameters (e.g., damping, step control, restart policies) and on the sensitivity of the iteration map to perturbations. This sensitivity is especially pronounced in high-order and parameterized iterative schemes, where the convergence dynamics may change significantly with the choice of internal parameters~\cite{petkovic2014multipoint}. Consequently, there is growing interest in diagnostic layers that can (i) map stability and reliability across parameter domains, and (ii) provide early predictions of solver behavior to reduce wasted computation and enable adaptive control.

In many biomedical and physiological modelling tasks, nonlinear equation solvers are embedded within larger computational pipelines, for example in parameter estimation for pharmacokinetic models~\cite{Kumar2024,Yan2014}, inverse problems in bioheat transfer, electrophysiological simulations, or models of neural dynamics~\cite{keener2009mathematical,quarteroni2007numerical}. In these settings, the reliability of the numerical solver plays a critical role: unstable iterations or convergence to non-physical solutions may propagate through the modelling workflow and lead to misleading quantitative conclusions~\cite{Creswell2024}. Moreover, biomedical models are often highly sensitive to parameter variations, making the behavior of the solver difficult to predict a priori.  Biological models spanning from biochemical reaction networks to disease spread simulations, often exhibit complex behaviors where small input perturbations result in substantial changes to outputs. This sensitivity stems from the inherent complexity of biological systems, which are often modeled with large numbers of free parameters that significantly affect model behavior. There are three key reasons for the high sensitivity to parameters as follows. (i) The nonlinear dynamics: biological processes, such as metabolic pathways or gene regulation, often involve nonlinear relationships (e.g., sigmoid curves), meaning small changes in a parameter can lead to a shift from negligible to maximum impact. (ii) The large number of connections between network nodes meaning that small changes in one parameter can cascade through networks, amplifying over time and space, causing significant phenotypic or ecological shifts. (iii) The structural sensitivity meaning that biological models are often highly sensitive to the specific mathematical functions used to describe processes (like growth or mortality), not just the parameter values, which can lead to divergent predictions.
These challenges motivate the development of diagnostic layers that provide transparent and interpretable indicators of solver reliability. By detecting instability early and guiding parameter selection, such tools can contribute to more trustworthy computational pipelines in biomedical modelling.

Machine learning (ML) has recently become a prominent tool for enhancing nonlinear solvers and numerical algorithms more broadly~\cite{karniadakis2021physics}. The literature largely clusters into four complementary directions. First, learned warm-starts use regression models (often neural networks) to provide high-quality initial guesses and improve basin selection, which can significantly reduce iterations and failure rates in simultaneous or high-order schemes \cite{shams2023inverse,shams2024qanalogues,freitas2021nnroots,zuo2025coldstart}. Second, surrogate-assisted solving replaces expensive solver components with learned approximations, for example through physics-informed neural networks (PINNs) or operator-learning frameworks \cite{raissi2019physics,lu2021deeponet}. Third, solver control and adaptivity uses machine learning to tune discretization parameters or solver hyperparameters during runtime \cite{hutter2019automated}. Fourth, early convergence or failure prediction analyzes short segments of iteration traces to forecast solver behavior and trigger restarts or parameter adjustments, often relying on indicators derived from the trajectory of the underlying iteration dynamics.

A closely related line of work explores trajectory-based stability indicators that characterize the dynamical behavior of iterative processes and can serve as informative signals for data-driven or machine-learning-based diagnostics. In particular, Lyapunov-exponent-based diagnostics and finite-time stability measures have long been used to characterize transient instability and chaotic behavior in nonlinear dynamical systems \cite{wolf1985determining}. More recently, machine learning approaches have been proposed to estimate Lyapunov exponents directly from time-series data \cite{velichko2025lle}. Related trajectory-driven diagnostics based on micro-ensemble perturbation analysis have also been studied for analyzing local stability and synchronization phenomena in nonlinear systems \cite{velichko2025microensemble}.

In the context of iterative numerical algorithms, these ideas naturally motivate trajectory-based stability profiling approaches. Iterative solvers generate sequences of approximations whose evolution can be interpreted as a discrete dynamical system, making them particularly amenable to trajectory-based diagnostics. In particular, our recent study \cite{shams2026contractivity} introduced a finite-time contractivity profiling framework based on a kNN micro-series estimator of a largest Lyapunov exponent (LLE) proxy, enabling systematic exploration of parameter spaces in parallel root-finding schemes. That work demonstrated that profile-derived metrics such as $S_{\min}$ and $S_{\mathrm{mom}}$ can be used to construct stability landscapes over $(\alpha,\beta)$ parameter grids and to identify regions associated with faster and more reliable convergence. Similar trajectory-driven diagnostics have also been explored for parameter optimization in parallel iterative methods through finite-time contraction profiling \cite{shams2026steplog}.

While these approaches provide systematic tools for parameter optimization and stability analysis, they do not address a complementary question that arises in practical computational workflows: whether the reliability of a solver configuration can be predicted early during the iteration process. In contrast to profiling-based approaches that require the full trajectory to assess solver stability, the present work investigates whether reliability can be predicted early from short prefixes of trajectory-derived stability indicators, thereby enabling lightweight AI-assisted diagnostics during solver execution. Importantly, the indicators used in this study retain a clear dynamical interpretation, providing interpretable signals linking solver behavior to stability properties of the underlying iteration process.

Despite these advances, several broader challenges remain in the integration of machine learning with numerical solvers. Generalization beyond the training distribution is often underexplored, interpretability may be limited, and empirical gains can sometimes be overstated when evaluated against weak baselines or under selective reporting~\cite{mcgreivy2024weak}. These issues motivate approaches that retain a clear physical or dynamical interpretation and provide transparent diagnostics rather than black-box prescriptions.

The present paper addresses this gap by introducing an interpretable AI-augmented reliability diagnostic framework for a two-parameter parallel root-finding scheme. The method builds on kNN--LLE proxy profiles derived from early iteration dynamics and formulates a multi-horizon learning problem: can a profile-based reliability score $S_{\mathrm{mom}}$ be predicted accurately from only the first few raw proxy values? The answer is affirmative. Across both random and structured center-to-periphery splits over the $(\alpha,\beta)$ grid, we show that strong predictive performance is achieved from surprisingly short prefixes, with prediction quality improving rapidly as the horizon increases and saturating thereafter. Because the learning layer operates on interpretable stability indicators rather than raw trajectories, the resulting framework remains closely connected to the dynamical analysis of solver behavior.

The contributions of this work are fourfold:
\begin{enumerate}
\item We develop a unified computational pipeline that couples kNN--LLE proxy profiling, profile-derived reliability scoring, and multi-horizon supervised learning over a two-parameter solver domain.
\item We show that early prefixes of the raw proxy profile already contain substantial information about the final reliability score $S_{\mathrm{mom}}$, enabling prediction before the characteristic minimum of the smoothed profile is observed.
\item We compare multiple standard regression models under both random and center-to-periphery splits, quantifying not only in-distribution accuracy but also extrapolative generalization across the parameter landscape.
\item We provide interpretable visual and statistical validation through profile examples, timing histograms, metric curves, and test-only heatmap reconstructions, thereby clarifying how the learned diagnostic layer relates to the underlying solver dynamics.
\end{enumerate}

The remainder of the paper is organized as follows. Section~\ref{sec:methodology} introduces the proposed framework, including the underlying iterative scheme, the construction of the kNN--LLE proxy profiles, the profile-based reliability metrics, and the multi-horizon learning protocol. Section~\ref{sec:experimental_setup} presents the experimental design, including the parameter grid, profiling hyperparameters, and evaluation procedure. Section~\ref{sec:results} reports the main empirical findings, combining profile-level illustrations, timing statistics, multi-horizon regression results, and heatmap-level validation under both random and center-to-periphery splits. Section~\ref{sec:discussion} interprets these results from the viewpoints of early diagnosis, computational savings, interpretability, and generalization. Finally, Section~\ref{sec:conclusion} summarizes the main conclusions and outlines directions for future work.

\section{Methodology}
\label{sec:methodology}

\begin{figure}[t]
\centering
\resizebox{\linewidth}{!}{%
\begin{tikzpicture}[
  font=\small,
  box/.style={
    rounded corners=2mm,
    draw=black!45,
    very thick,
    align=center,
    minimum height=10mm,
    inner sep=2mm,
    text width=4.5cm,
    fill=#1,
    blur shadow={shadow xshift=0.8mm, shadow yshift=-0.8mm, shadow blur steps=4}
  },
  arr/.style={-{Latex[length=3mm,width=2mm]}, very thick, draw=black!70}
]

\node[box=blue!14]  (A) {\textbf{Parameter grid}\\\& solver ensemble runs};
\node[box=cyan!14,  right=7mm of A] (B) {\textbf{kNN--LLE proxy profile}\\$\lambda_1(t)$ (raw, smoothed)};
\node[box=green!16, right=7mm of B] (C) {\textbf{Profile metrics}\\$S_{\min},\,S_{\mathrm{mom}},\,t_{\min},\ldots$};

\node[box=orange!18, below=11mm of B] (D) {\textbf{ML dataset}\\multi-horizon prefixes $X_T$};
\node[box=red!14,    right=7mm of D] (E) {\textbf{Early reliability}\\prediction $\widehat{S}_{\mathrm{mom}}(T)$};
\node[box=purple!14, right=7mm of E] (F) {\textbf{Validation}\\metric curves \& test-only heatmaps};

\draw[arr] (A) -- (B);
\draw[arr] (B) -- (C);
\draw[arr] (B) -- (D);
\draw[arr] (D) -- (E);
\draw[arr] (E) -- (F);

\end{tikzpicture}%
}
\caption{Conceptual pipeline of the proposed reliability diagnostic. Solver ensembles are generated on a $(\alpha,\beta)$ grid, yielding kNN--LLE proxy profiles and interpretable profile metrics. Multi-horizon prefixes form an ML dataset used for early prediction of $S_{\mathrm{mom}}(T)$, followed by validation via metric curves and test-only heatmaps.}
\label{fig:pipeline}
\end{figure}
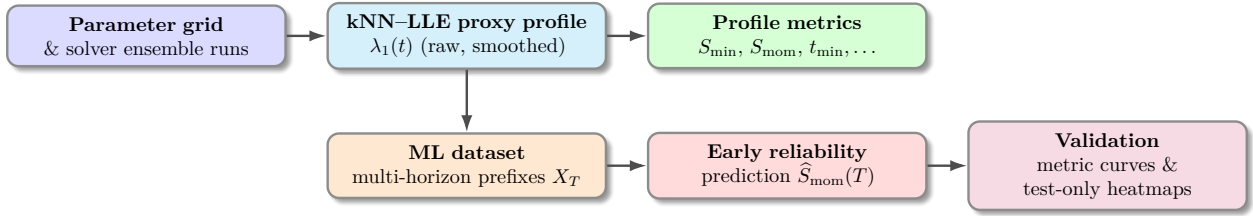

Figure~\ref{fig:pipeline} summarizes the overall computational pipeline used in this work. The workflow begins with profile generation on a two-parameter grid, where repeated solver runs are used to construct kNN--LLE proxy profiles and profile-derived reliability metrics. These profile traces are then converted into multi-horizon datasets by truncating the raw proxy profile at increasing prefix lengths $T$. Standard regression models are trained to predict the final reliability score $S_{\mathrm{mom}}$ from these early prefixes, enabling quantitative analysis of how soon solver reliability becomes predictable from early dynamics. The learned models are finally evaluated through metric curves and test-only heatmap reconstructions over the parameter domain.

\subsection{Problem Setting and Iterative Scheme (Brief)}
We consider nonlinear root-finding problems of the form
\begin{equation}
f(z)=0,
\end{equation}
where $f:\mathbb{R}\to\mathbb{R}$ (and, more generally, $f:\mathbb{R}^d\to\mathbb{R}^d$) may be nonlinear and parameter-dependent.
The computational backbone is a derivative-free, parallel, two-parameter iterative scheme whose behavior depends on two control parameters $(\alpha,\beta)$.
For each fixed $(\alpha,\beta)$, the solver generates an iterate sequence $\{z_k\}_{k\ge 0}$ from an initialization $z_0$.
The goal of the present study is to build an interpretable reliability layer around this solver rather than replacing the numerical scheme.
Auxiliary solver-level validation under varied initial guesses, together with a comparison against existing parallel schemes, is reported in Appendix~\ref{app:solver_validation}. Those appendix experiments are carried out at a single representative parameter pair, $(\alpha,\beta)=(-0.1,4.0)$, selected from the stable region of the $(\alpha,\beta)$ parameter map.

\subsection{kNN--LLE Proxy Profile Construction}
For each solver run we define a scalar micro-series that summarizes early iterative dynamics. In this study we use the log step norm
\begin{equation}
y_k=\log\bigl(\|z_{k+1}-z_k\|\bigr), \qquad k=0,1,\ldots,K-1,
\end{equation}
computed over a fixed horizon of $K$ iterations. For a given $(\alpha,\beta)$, an ensemble of independent runs yields a collection of micro-series.

To obtain a lightweight proxy for local stability, we estimate a kNN-based largest Lyapunov exponent (LLE) indicator from short windows of the micro-series.
For each window end index $t_{\mathrm{end}}$, delay-embedded vectors are built using a look-back length $L$ and multiple prediction horizons $h\in\{h_{\min},\ldots,h_{\max}\}$.
The slope of the multi-horizon forecast error acts as a stability proxy, producing a profile
\begin{equation}
\lambda_1(t_{\mathrm{end}}), \qquad t_{\mathrm{end}}=1,2,\ldots,W,
\end{equation}
where negative values indicate locally contractive dynamics and positive values indicate locally expansive behavior.
A trailing-mean filter with window size $W_{\mathrm{smooth}}$ is applied to obtain a smoothed profile $\tilde{\lambda}_1(t)$.

\subsection{Profile-Based Reliability Metrics}
From $\tilde{\lambda}_1(t)$ we compute interpretable summary metrics capturing both depth and extent of contractive behavior. Let
\begin{equation}
t_{\min}=\arg\min_{t\in[t_{\mathrm{start}},t_{\mathrm{end}}]}\tilde{\lambda}_1(t), \qquad
y_{\min}=\min_{t\in[t_{\mathrm{start}},t_{\mathrm{end}}]}\tilde{\lambda}_1(t).
\end{equation}
We define
\begin{equation}
S_{\min}=\frac{-y_{\min}}{t_{\min}+\varepsilon},
\end{equation}
and a negative-mass moment score
\begin{equation}
M_0=\sum_t \max\{0,-\tilde{\lambda}_1(t)\}, \quad
\bar{t}=\frac{\sum_t t\,\max\{0,-\tilde{\lambda}_1(t)\}}{M_0+\varepsilon}, \quad
S_{\mathrm{mom}}=\frac{M_0}{\bar{t}+\varepsilon},
\end{equation}
where sums are taken over a fixed evaluation interval and $\varepsilon>0$ is a small numerical constant.
Intuitively, larger $S_{\mathrm{mom}}$ indicates earlier and stronger contractivity.

\subsection{Dataset Construction for Multi-Horizon Learning}
To predict $S_{\mathrm{mom}}$ without running the profiling procedure until the full minimum is observed, we define multi-horizon learning tasks using early prefixes of the raw profile.
For each grid point $(\alpha,\beta)$ and each horizon $T$ we define
\begin{equation}
X_T=[\lambda_1^{\mathrm{raw}}(1),\lambda_1^{\mathrm{raw}}(2),\ldots,\lambda_1^{\mathrm{raw}}(T)],
\end{equation}
with target $y=S_{\mathrm{mom}}$ computed from the smoothed profile. The set of horizons is $\mathcal{T}=\{1,1+\Delta_T,1+2\Delta_T,\ldots\}$ up to the available profile length.

\subsection{Train/Test Splits: Random vs.\ Center-to-Periphery Generalization}
We evaluate generalization using two complementary splitting strategies.

\paragraph{Random split.}
All grid points are shuffled and divided into training and test sets according to a fixed test fraction. This split assesses average-case generalization when training and test samples follow the same parameter distribution.

\paragraph{Center-to-periphery split.}
To probe extrapolation across parameter regimes, we compute each point's normalized squared distance to the grid center in $(\alpha,\beta)$ and allocate the most central fraction to training,
with peripheral points used for testing. This split stresses out-of-distribution behavior near regime boundaries.

\subsection{Learning Models and Evaluation (Brief)}
For each horizon $T$, we train a regressor to predict the reliability score $S_{\mathrm{mom}}$ from the early-prefix feature vector $X_T$. To cover complementary inductive biases and complexity levels, we consider five standard model families: (i) k-nearest neighbors (kNN) regression, which provides a local, instance-based predictor well suited to regime-dependent behaviors; (ii) Ridge regression (linear model with $\ell_2$ regularization), which offers a stable baseline under correlated features; (iii) Elastic Net (combined $\ell_1$/$\ell_2$ regularization), which balances stability with sparse feature selection; (iv) Random Forest regression, an ensemble of decision trees capturing nonlinear interactions and heterogeneous regimes; and (v) Gradient Boosting regression, a sequential tree-ensemble that often achieves high accuracy by fitting additive corrections. For each $T$, models are trained independently to reflect the multi-horizon (early-exit) setting and to quantify how predictive quality evolves as more early profile information becomes available.

Model performance is reported using standard regression criteria on held-out test data. In the final manuscript we focus on three complementary metrics: MAE (mean absolute error) as a robust measure of typical absolute deviation, RMSE (root mean squared error) as a penalty for larger errors, and $R^2$ (coefficient of determination) as a scale-free indicator of explained variance. These metrics together characterize both absolute accuracy and overall predictive fidelity across the parameter grid.

\section{Experimental Setup}
\label{sec:experimental_setup}

\subsection{Parameter Grid and Computational Budget}
All experiments are performed over a two-parameter control grid $(\alpha,\beta)$, with
\begin{equation}
\alpha \in [-3,5], \qquad \beta \in [-2,4],
\end{equation}
sampled on a uniform $60\times 60$ grid. For each grid point, we generate an ensemble of $N_{\mathrm{runs}}=1000$ independent solver trajectories (distinct initializations) and run the iterative scheme for $K=200$ iterations per trajectory. Unless stated otherwise, a fixed random seed is used to ensure exact reproducibility across runs.

\subsection{Micro-Series Definition and Stabilization Settings}
From each trajectory we construct the micro-series
\begin{equation}
y_k=\log\bigl(\|z_{k+1}-z_k\|\bigr), \qquad k=0,\ldots,K-1,
\end{equation}
which is used as the input signal for kNN--LLE proxy estimation. To avoid numerical artifacts in extremely small step regimes, a log-domain tail floor is applied (freezing the tail after the floor is reached). In addition, an adaptive stabilization mechanism is enabled only to prevent explosive divergence (e.g., hard bounds for divergence markers and a step cap when stabilization is active), while preserving the raw dynamics in non-explosive regimes.

\subsection{kNN--LLE Proxy Hyperparameters}
The kNN--LLE proxy profile $\lambda_1(t_{\mathrm{end}})$ is computed using the following fixed hyperparameter configuration:
\begin{itemize}
    \item look-back length (embedding depth): $L=5$,
    \item forecast horizons: $h \in \{1,2,3,4,5\}$ (i.e., $h_{\min}=1$, $h_{\max}=5$, step $1$),
    \item kNN neighborhood size: $k=3$,
    \item internal train/test split for the micro-series predictor: $0.60/0.40$.
\end{itemize}
The raw profile is smoothed by a trailing-mean filter with window size $W_{\mathrm{smooth}}=4$ to reduce small oscillations.

\subsection{Metric Window and Reliability Targets}
Profile-based metrics are evaluated on a fixed index window of the proxy profile. Specifically, we consider
\begin{equation}
t \in [t_{\mathrm{start}},t_{\mathrm{end}}]=[10,200],
\end{equation}
where the effective start index is constrained by the minimum micro-series length required by the embedding and multi-horizon settings, i.e.,
\begin{equation}
t_{\mathrm{start}}=\max\{10,\;L+h_{\max}\},
\end{equation}
with look-back length $L=5$ and $h_{\max}=5$ in the present setup. The raw proxy profile $\lambda_1^{\mathrm{raw}}(t)$ is additionally smoothed using a trailing-mean filter of window size $W_{\mathrm{smooth}}=4$, yielding $\tilde{\lambda}_1(t)$, which is used for metric computation. The primary learning target is the smoothed-profile reliability score $S_{\mathrm{mom}}$ (defined in the Profile-Based Reliability Metrics subsection), interpreted as a scalar indicator of early and strong contractive behavior and used to rank $(\alpha,\beta)$ configurations in terms of expected solver reliability.

\subsection{Multi-Horizon Learning Protocol}
For each grid point, learning features are constructed from the raw profile prefix $X_T=[\lambda^{\mathrm{raw}}_1(1),\ldots,\lambda^{\mathrm{raw}}_1(T)]$ with horizons scanned as
\begin{equation}
T \in \{1, 1+\Delta_T, 1+2\Delta_T,\ldots\},
\end{equation}
up to the available profile length (with a fixed step $\Delta_T$). Each horizon $T$ defines an independent regression task mapping $X_T \mapsto S_{\mathrm{mom}}$, enabling an "early-exit" assessment of predictive performance as a function of the amount of early information used.

\subsection{Train/Test Splits and Evaluation}
We evaluate generalization using two complementary train/test splitting strategies: (i) a random split over all grid points with a fixed test fraction of $40\%$, and (ii) a center-to-periphery split, where training uses the most central fraction of grid points in $(\alpha,\beta)$ space and testing uses peripheral points. For each split and each horizon $T$, model performance is reported on the held-out test set using MAE, RMSE, and $R^2$.

\section{Results}
\label{sec:results}

\subsection{Ground-Truth Distributions of $S_{\min}$ and $S_{\mathrm{mom}}$}
\begin{figure}[t]
\centering
\begin{minipage}{0.49\linewidth}
\centering
\includegraphics[width=\linewidth]{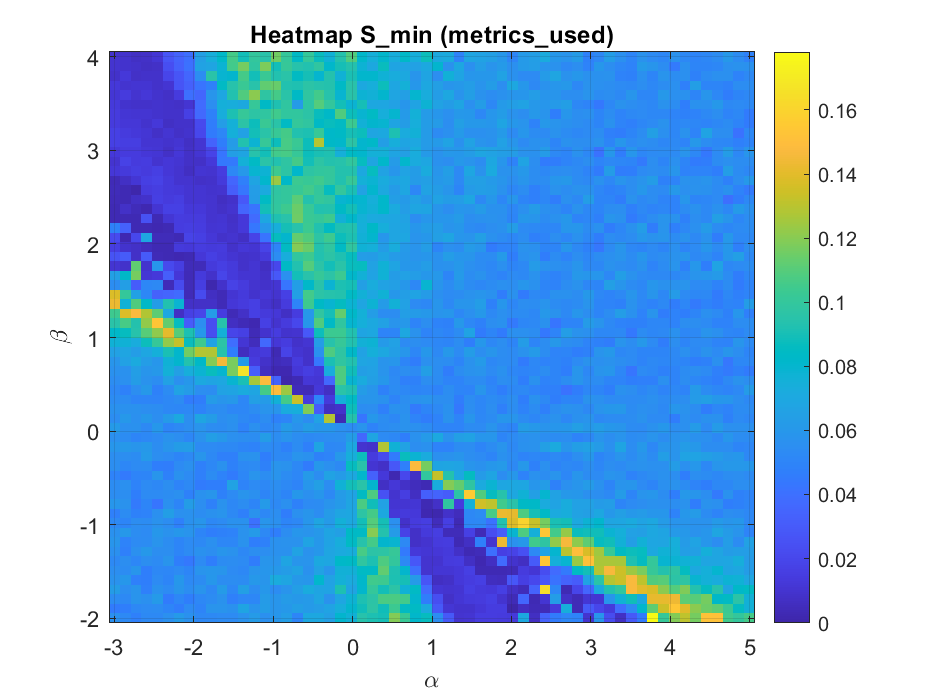}\\
(a) $S_{\min}$ heatmap.
\end{minipage}\hfill
\begin{minipage}{0.49\linewidth}
\centering
\includegraphics[width=\linewidth]{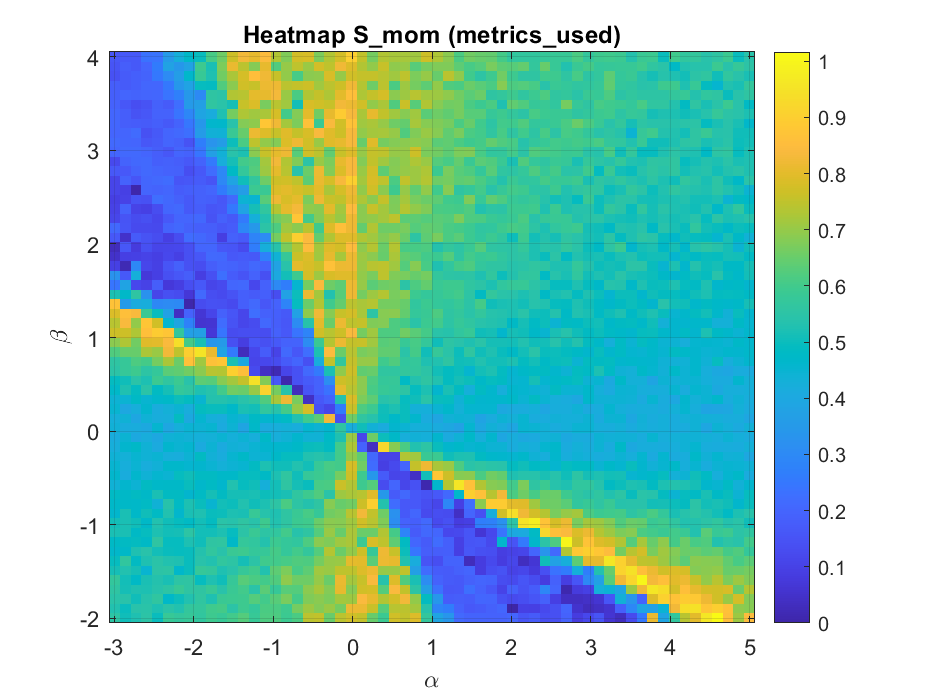}\\
(b) $S_{\mathrm{mom}}$ heatmap.
\end{minipage}
\caption{Ground-truth heatmaps of profile-based reliability metrics over the $(\alpha,\beta)$ grid: (a) $S_{\min}$ and (b) $S_{\mathrm{mom}}$.}
\label{fig:true_heatmaps}
\end{figure}

Figure~\ref{fig:true_heatmaps} reports the ground-truth heatmaps of the two profile-based reliability metrics over the $(\alpha,\beta)$ grid. The maps exhibit a highly structured landscape with distinct parameter regions associated with stronger contractive behavior (higher scores) versus weaker or delayed contractivity (lower scores). In particular, $S_{\mathrm{mom}}$ provides a clear, interpretable ranking of parameter configurations, highlighting a stable "good" region that is consistently separated from less reliable regimes.

\subsection{Representative kNN--LLE proxy profiles across the $S_{\mathrm{mom}}$ scale}
\begin{figure*}[t]
\centering

\begin{minipage}{0.32\linewidth}
\centering
\includegraphics[width=\linewidth]{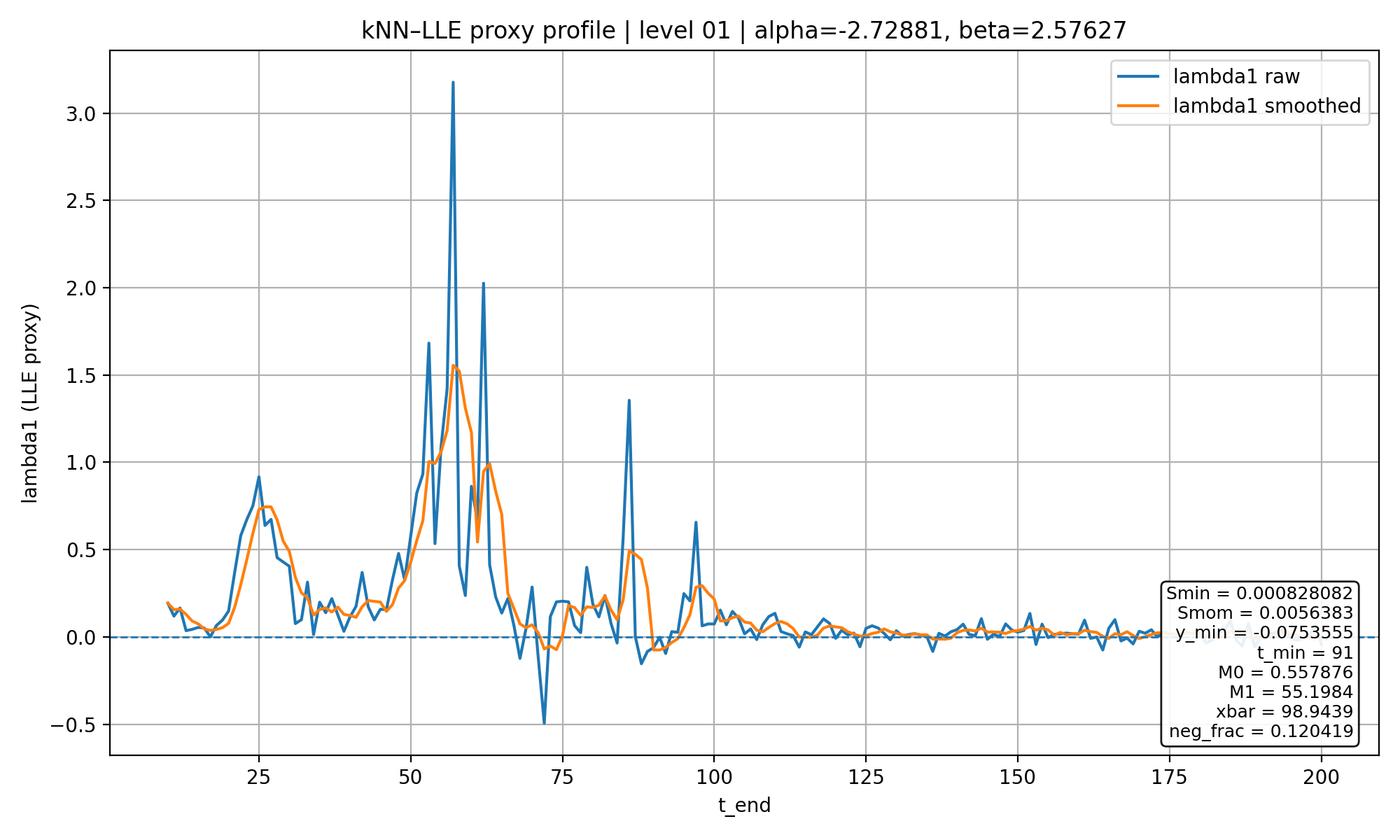}\\[-2pt]
(a) Level 01: $\alpha=-2.7288$, $\beta=2.5763$, $S_{\mathrm{mom}}=0.00564$.
\end{minipage}\hfill
\begin{minipage}{0.32\linewidth}
\centering
\includegraphics[width=\linewidth]{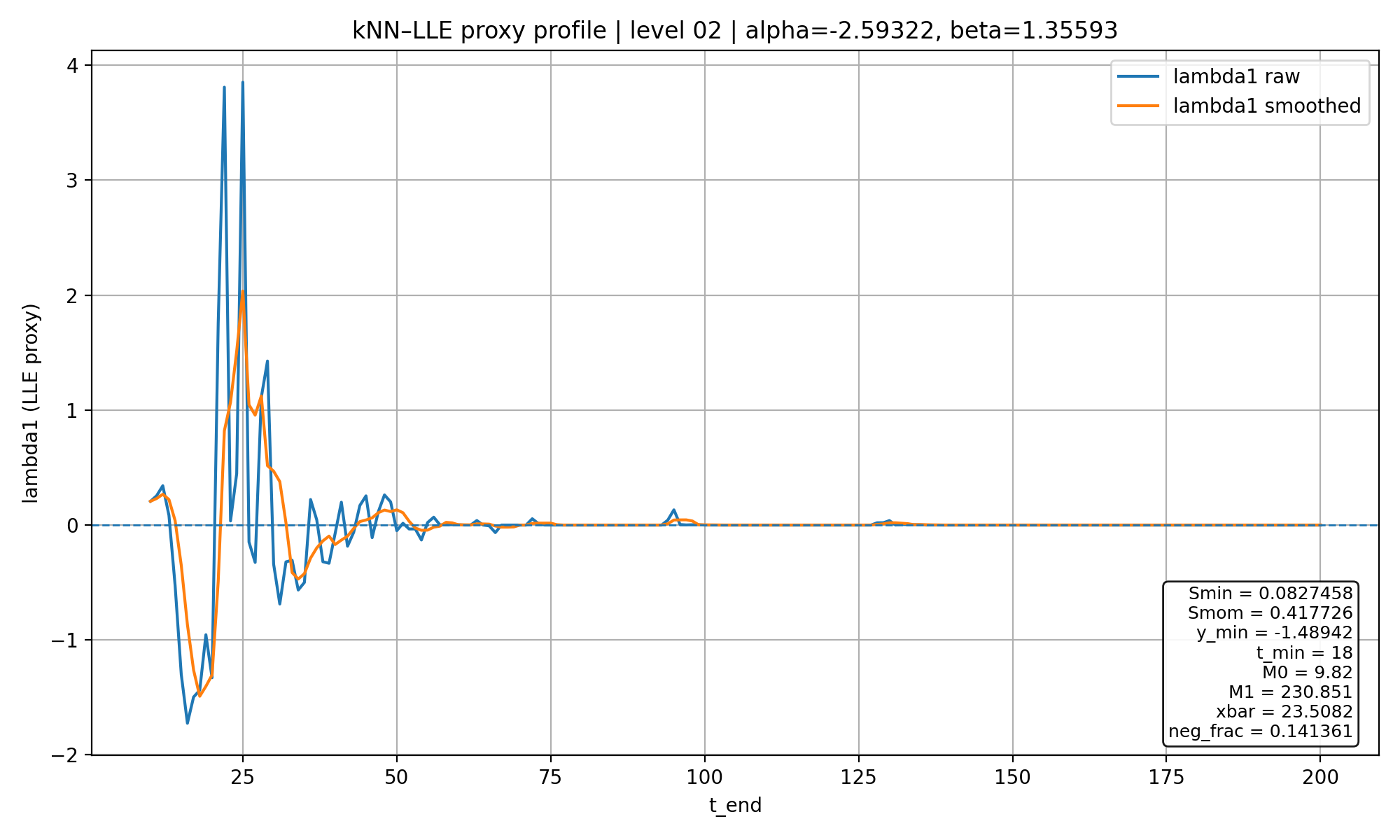}\\[-2pt]
(b) Level 02: $\alpha=-2.5932$, $\beta=1.3559$, $S_{\mathrm{mom}}=0.4177$.
\end{minipage}\hfill
\begin{minipage}{0.32\linewidth}
\centering
\includegraphics[width=\linewidth]{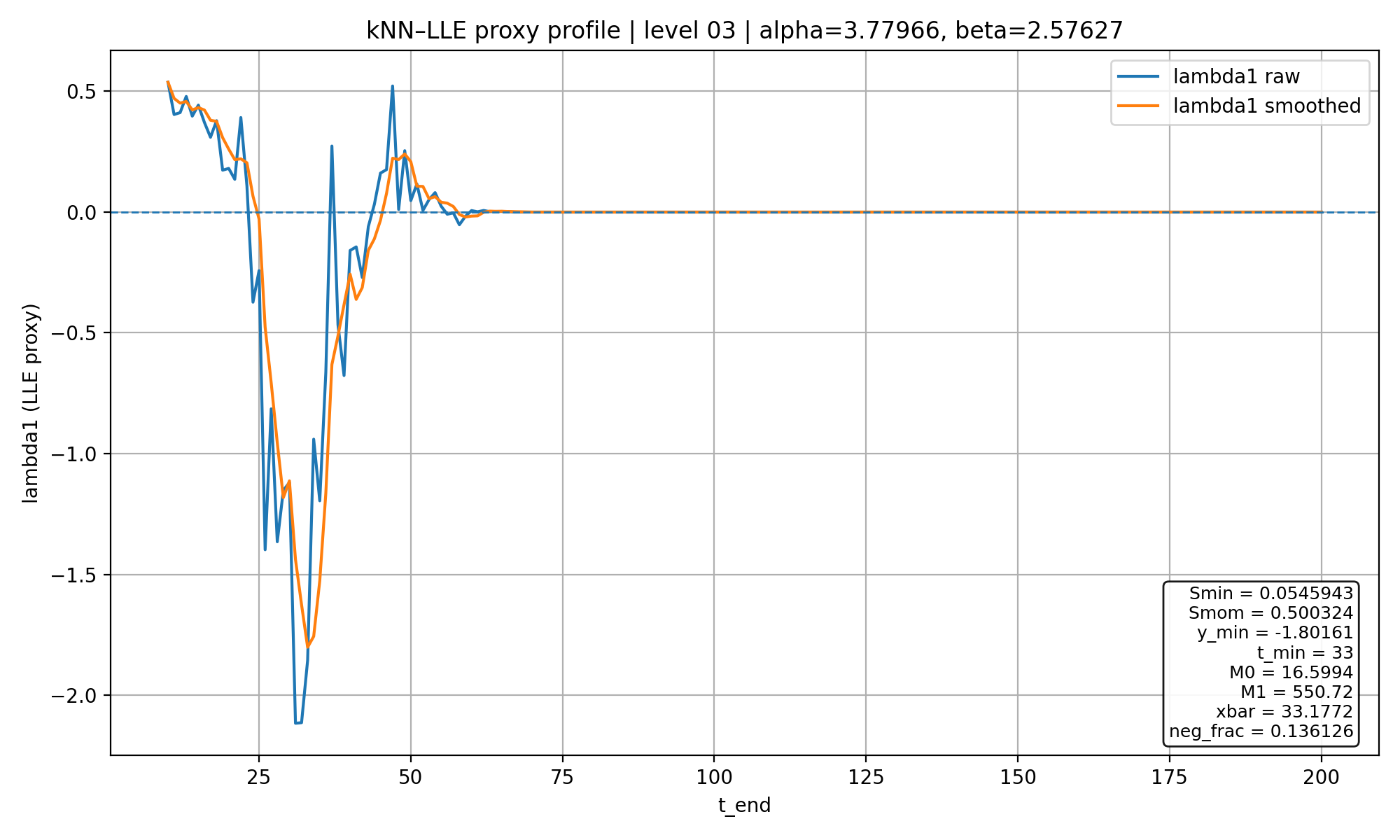}\\[-2pt]
(c) Level 03: $\alpha=3.7797$, $\beta=2.5763$, $S_{\mathrm{mom}}=0.5003$.
\end{minipage}

\vspace{6pt}

\begin{minipage}{0.32\linewidth}
\centering
\includegraphics[width=\linewidth]{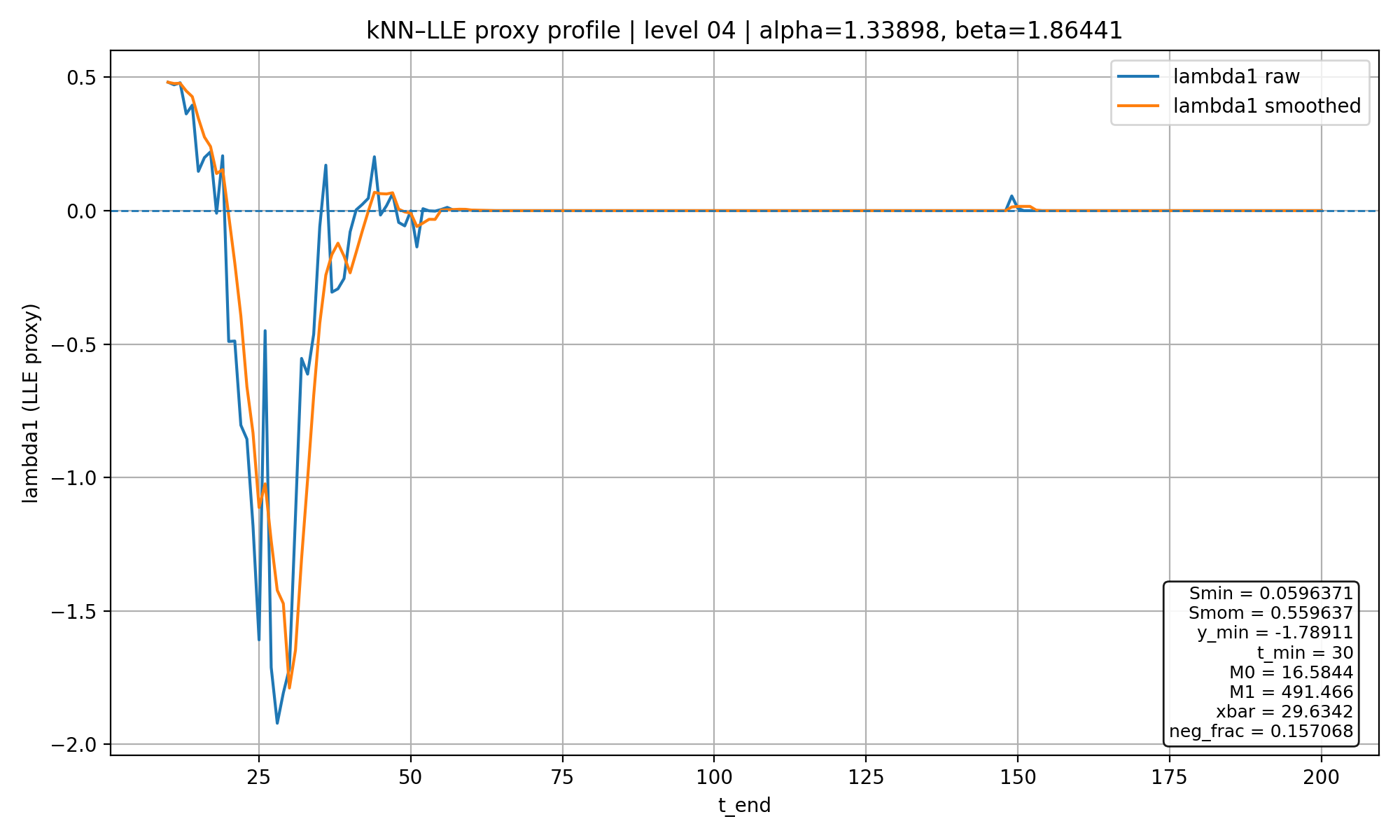}\\[-2pt]
(d) Level 04: $\alpha=1.3390$, $\beta=1.8644$, $S_{\mathrm{mom}}=0.5596$.
\end{minipage}\hfill
\begin{minipage}{0.32\linewidth}
\centering
\includegraphics[width=\linewidth]{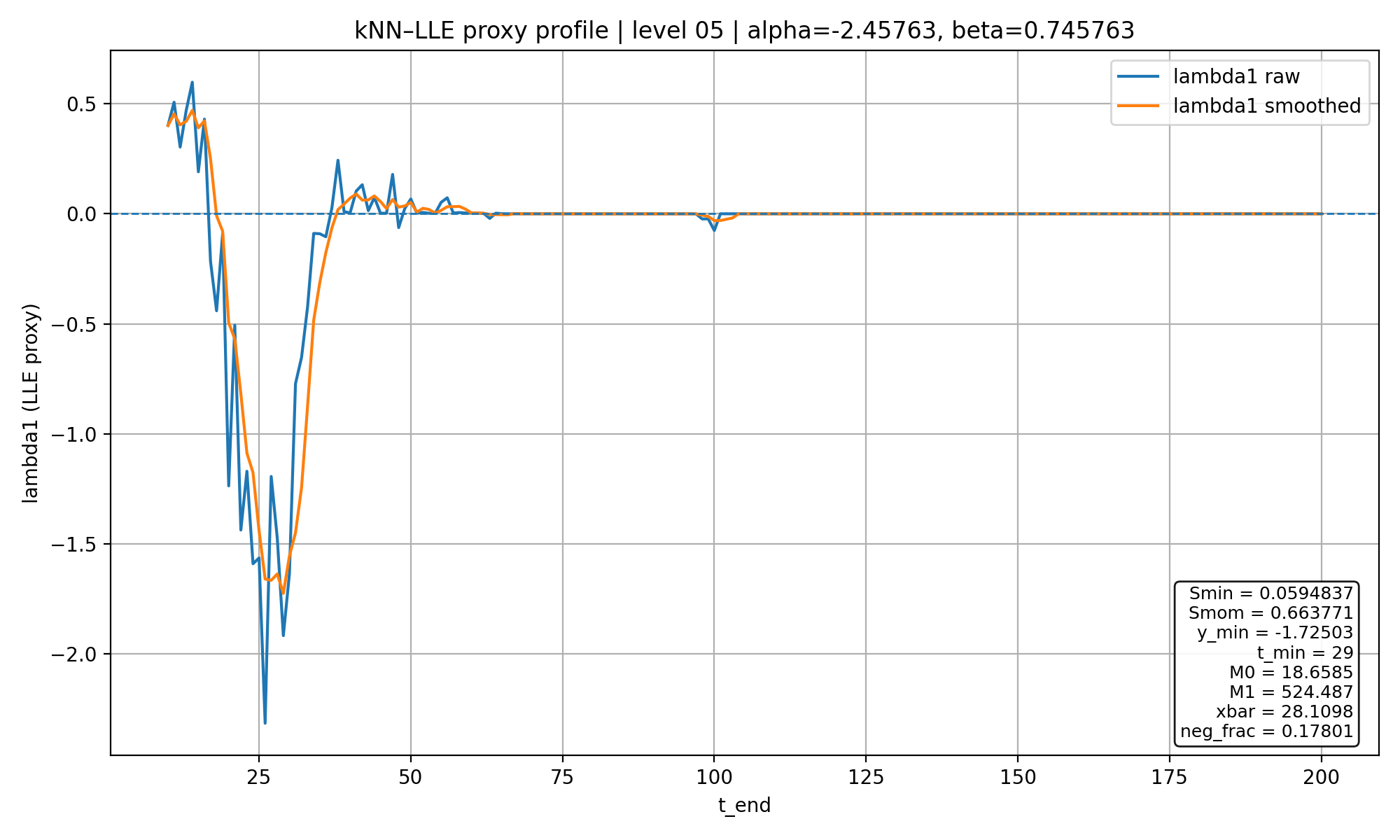}\\[-2pt]
(e) Level 05: $\alpha=-2.4576$, $\beta=0.7458$, $S_{\mathrm{mom}}=0.6638$.
\end{minipage}\hfill
\begin{minipage}{0.32\linewidth}
\centering
\includegraphics[width=\linewidth]{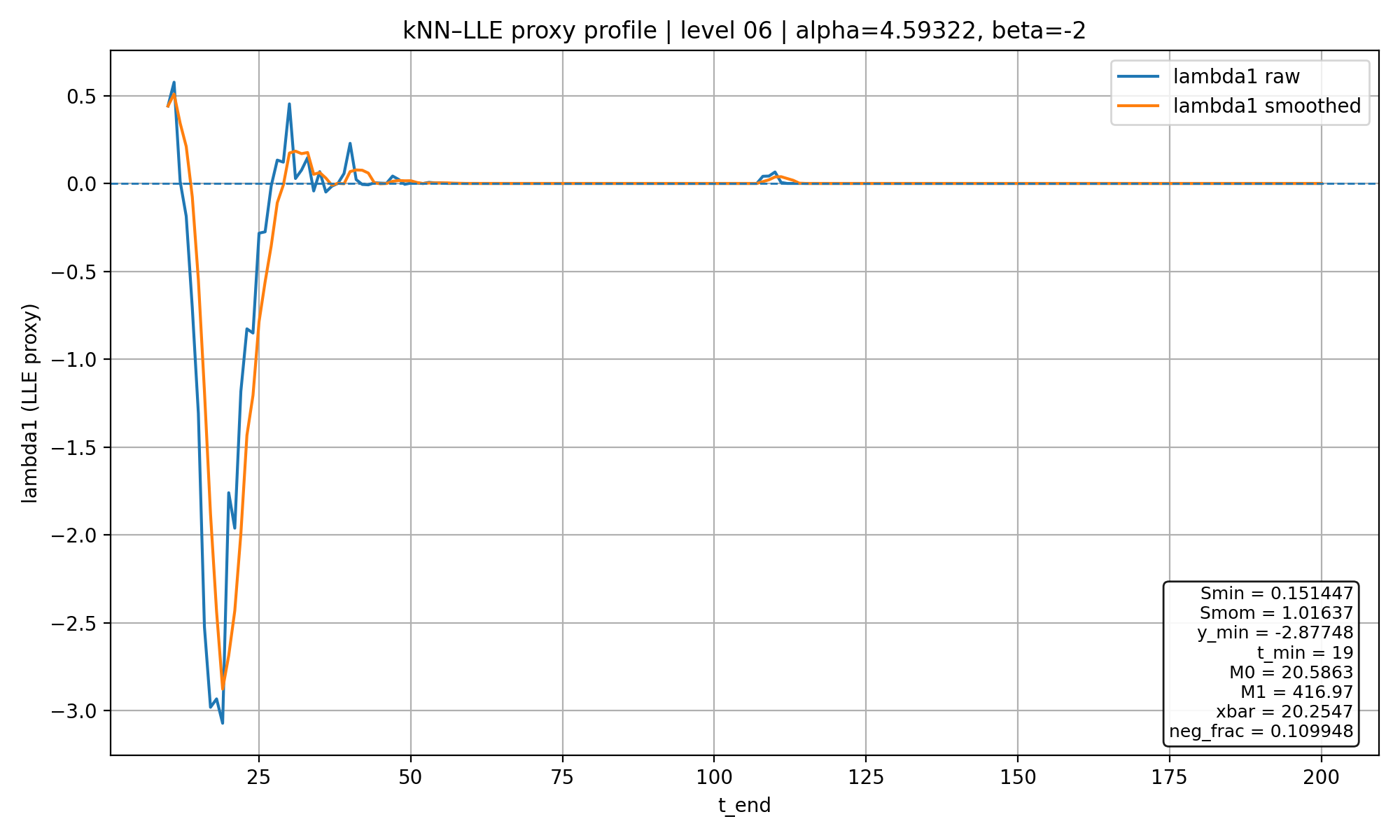}\\[-2pt]
(f) Level 06: $\alpha=4.5932$, $\beta=-2.0000$, $S_{\mathrm{mom}}=1.0164$.
\end{minipage}

\caption{Representative kNN--LLE proxy profiles (raw vs.\ smoothed) spanning increasing $S_{\mathrm{mom}}$ levels (a)--(f). Higher $S_{\mathrm{mom}}$ typically corresponds to earlier and deeper contractive dips (more negative values) in the smoothed profile, whereas low $S_{\mathrm{mom}}$ cases tend to exhibit weak or delayed contractivity.}
\label{fig:profiles_levels}
\end{figure*}

To qualitatively illustrate how the kNN--LLE proxy profile varies across the parameter domain, we select six representative grid points spanning the full observed range of $S_{\mathrm{mom}}$ (ordered by increasing $S_{\mathrm{mom}}$). Figure~\ref{fig:profiles_levels} shows raw and smoothed proxy profiles $\lambda_1(t_{\mathrm{end}})$ for each selected point.

A clear qualitative trend is observed: low-$S_{\mathrm{mom}}$ cases tend to exhibit predominantly non-contractive behavior (profiles largely nonnegative or only weakly negative), often with delayed and shallow minima (large $t_{\min}$). In contrast, high-$S_{\mathrm{mom}}$ cases display an earlier and deeper negative excursion, yielding a pronounced minimum and a larger negative mass concentrated at small $t_{\mathrm{end}}$. These visual patterns support the interpretability of $S_{\mathrm{mom}}$ as an "early-and-strong contractivity" indicator.

Importantly, the learning task addressed in this work is not to recompute $S_{\mathrm{mom}}$ from the full profile, but to predict $S_{\mathrm{mom}}$ from only an initial prefix of the raw curve (multi-horizon setting), i.e., before the minimum is necessarily observed.

\subsection{Timing Statistics of Contractivity Signatures}
To characterize when contractive behavior typically emerges, we analyze timing statistics derived from the smoothed proxy profile $\tilde{\lambda}_1(t)$. We define a "good" subset as the top $20\%$ of grid points ranked by $S_{\mathrm{mom}}$, yielding a threshold $S_{\mathrm{mom}}\ge 0.6638$ (720 points out of 3600). For this good subset, the median location of the profile minimum is $t_{\min}=21$ (corresponding to the profile index $T_{\min}\approx 12$ on the discrete $t_{\mathrm{end}}$ axis). This indicates that, in the good region, the characteristic contractivity minimum occurs relatively early, supporting the multi-horizon viewpoint where reliable discrimination can be achieved using only the initial portion of the profile.

\begin{figure}[t]
\centering
\begin{minipage}{0.49\linewidth}
\centering
\includegraphics[width=\linewidth]{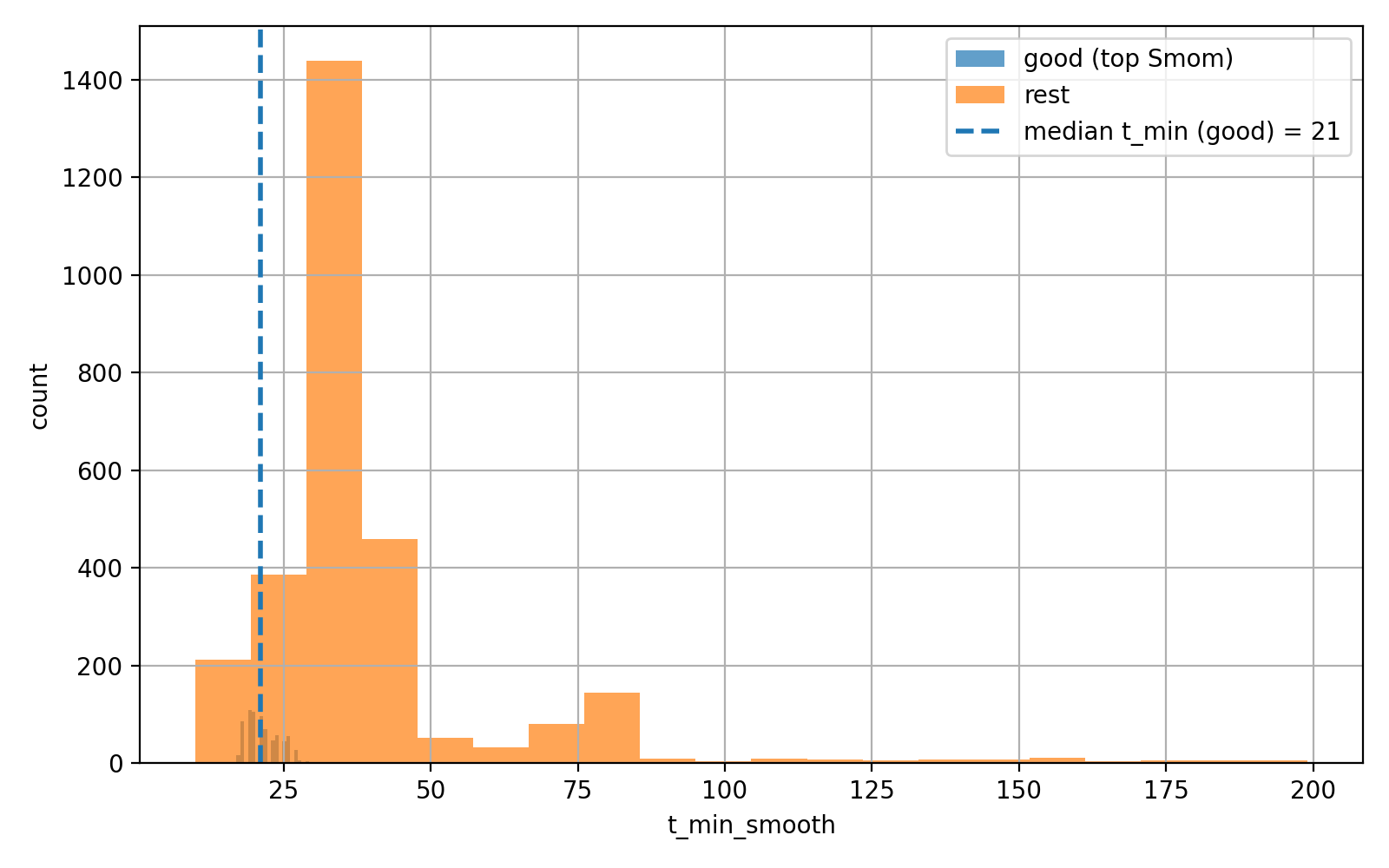}\\
(a) Histogram of $t_{\min}$ (smoothed profile).
\end{minipage}\hfill
\begin{minipage}{0.49\linewidth}
\centering
\includegraphics[width=\linewidth]{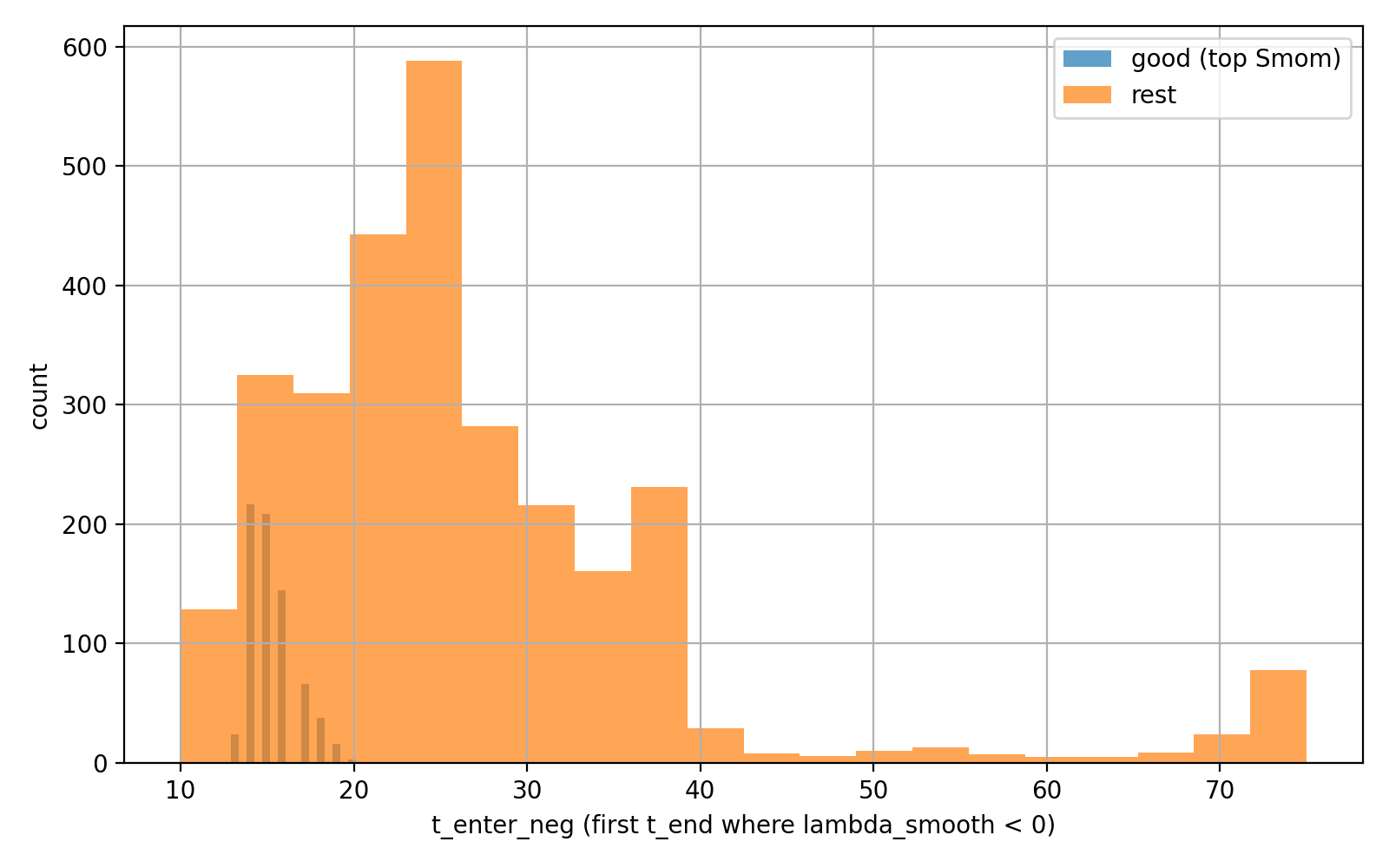}\\
(b) Histogram of $t_{\mathrm{enter\_neg}}$ (first $\tilde{\lambda}_1(t)<0$).
\end{minipage}
\caption{Timing statistics for the smoothed proxy profile $\tilde{\lambda}_1(t)$, comparing the good subset (top $20\%$ by $S_{\mathrm{mom}}$) against the rest: (a) distribution of the minimum location $t_{\min}$; (b) distribution of the first negative entry time $t_{\mathrm{enter\_neg}}$.}
\label{fig:timing_hists}
\end{figure}

Figure~\ref{fig:timing_hists} compares two timing distributions between the good subset and the remaining points: (i) the first time index at which $\tilde{\lambda}_1(t)$ becomes negative ($t_{\mathrm{enter\_neg}}$), and (ii) the minimum location $t_{\min}$. The good subset tends to enter the negative (contractive) regime earlier and exhibits a more concentrated distribution of $t_{\min}$, whereas the remaining region shows broader, often delayed timing patterns.

\subsection{Multi-Horizon Regression Results (Center-to-Periphery Split)}
We next evaluate the multi-horizon learning problem under the center-to-periphery split, where training samples are drawn from the most central region of the $(\alpha,\beta)$ grid and testing is performed on the periphery. This split probes out-of-distribution generalization across the parameter domain. All performance curves are reported as a function of the prefix length $T$ (number of early raw proxy values $\lambda^{\mathrm{raw}}_1(1{:}T)$ used as input), and are shown up to $T \le 3T_{\min}$, where $T_{\min}=12$ denotes the median minimum-location index of the good region (Figure~\ref{fig:timing_hists}).

\begin{figure}[t]
\centering
\begin{minipage}{0.32\linewidth}
\centering
\includegraphics[width=\linewidth]{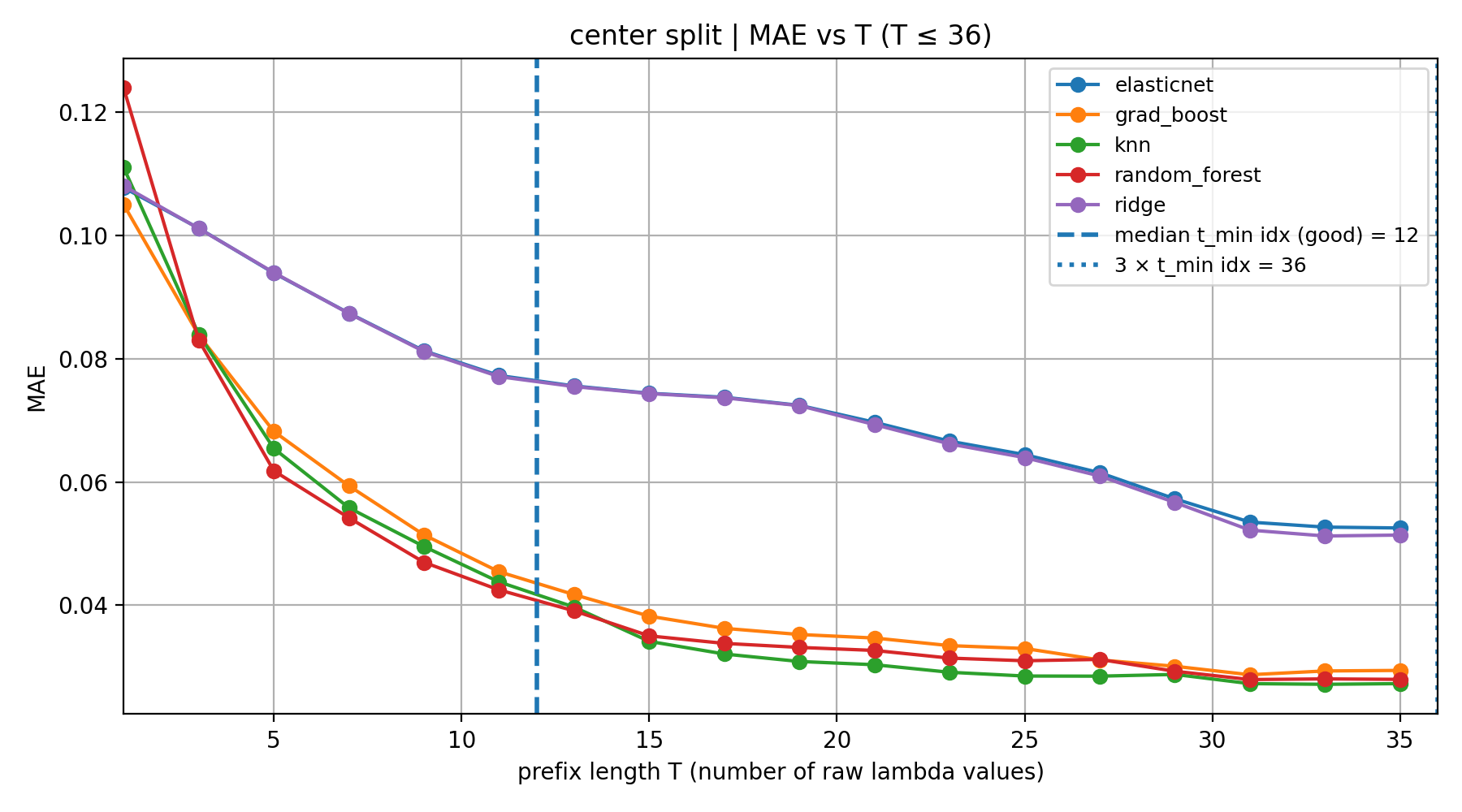}\\[-2pt]
(a) MAE vs.\ $T$.
\end{minipage}\hfill
\begin{minipage}{0.32\linewidth}
\centering
\includegraphics[width=\linewidth]{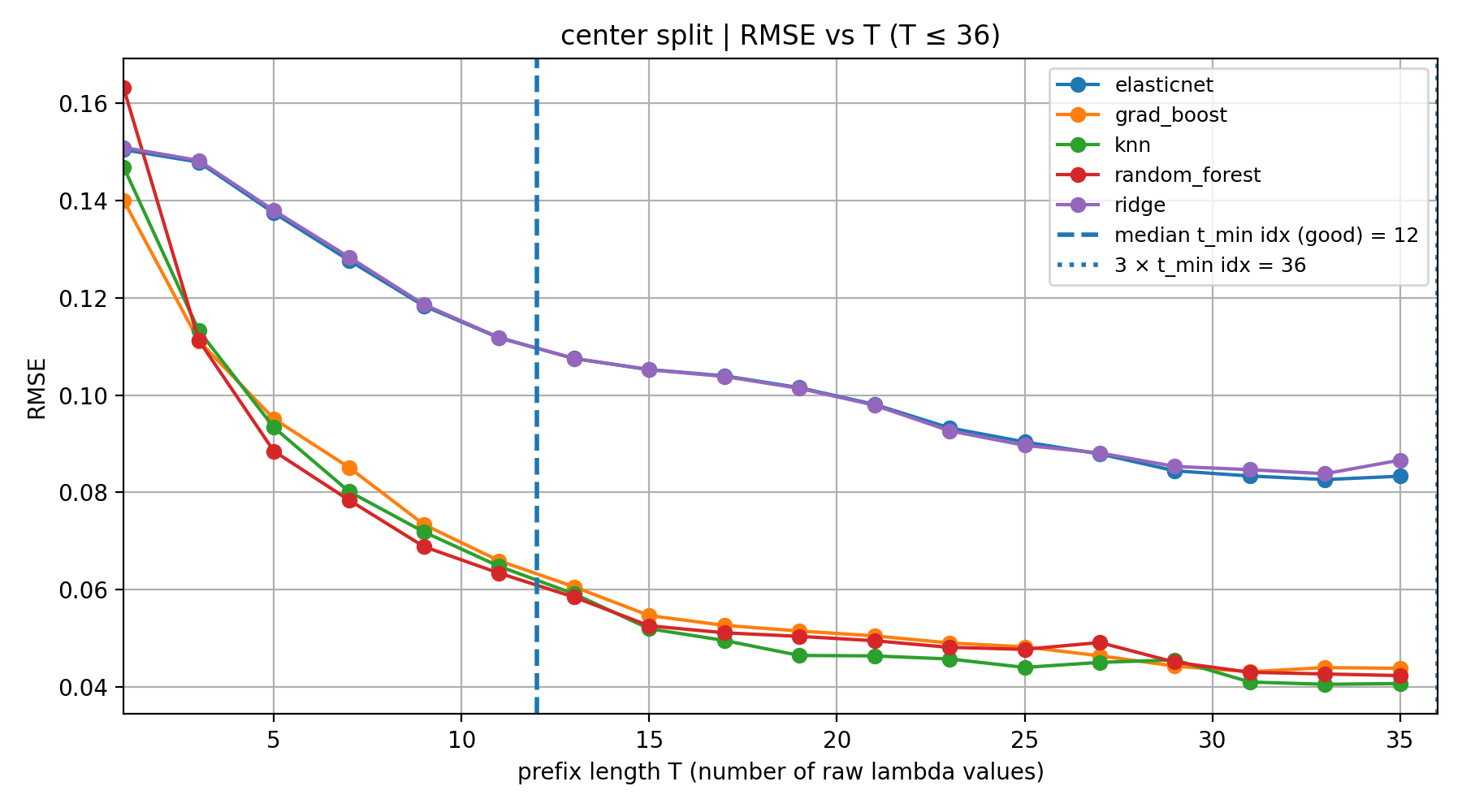}\\[-2pt]
(b) RMSE vs.\ $T$.
\end{minipage}\hfill
\begin{minipage}{0.32\linewidth}
\centering
\includegraphics[width=\linewidth]{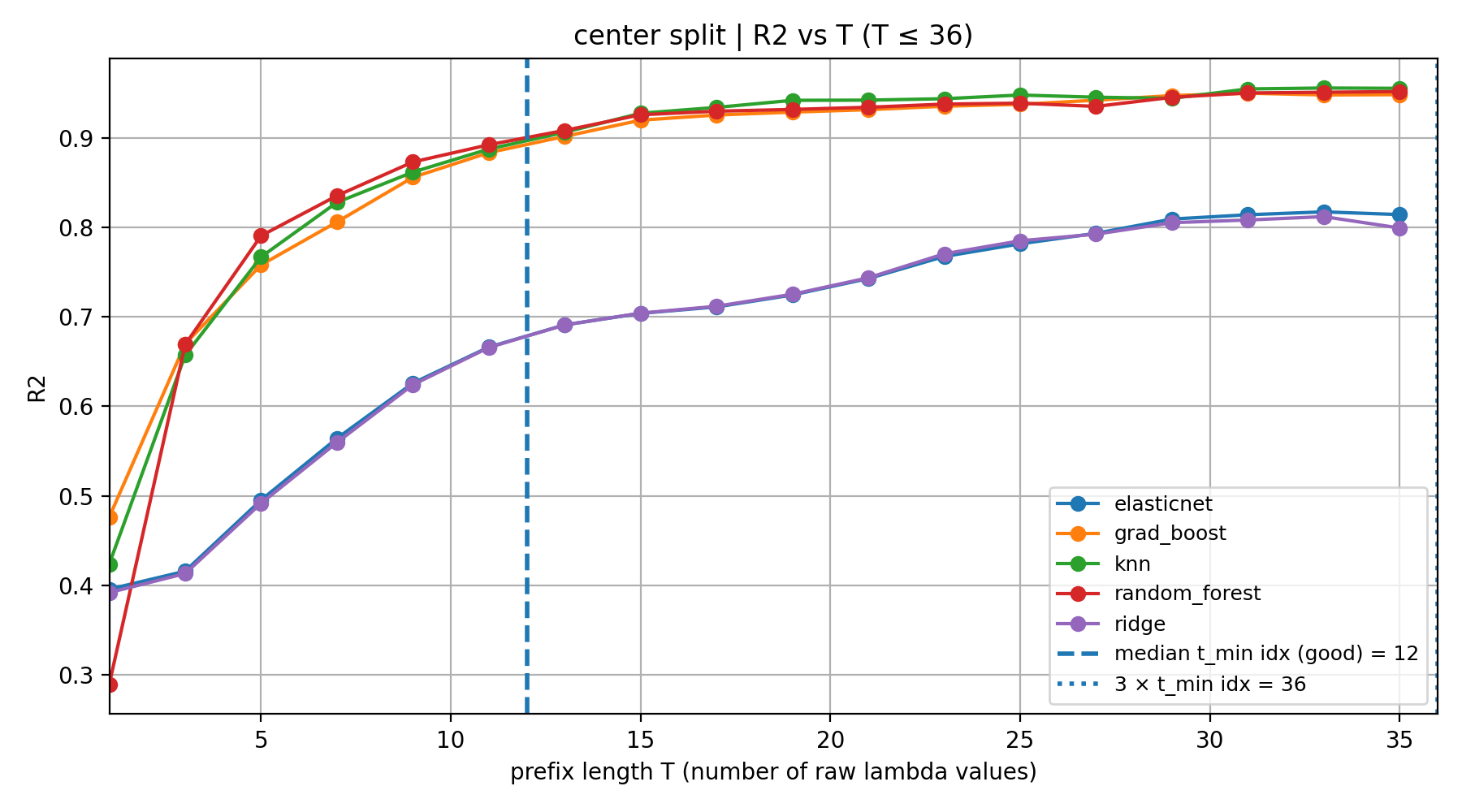}\\[-2pt]
(c) $R^2$ vs.\ $T$.
\end{minipage}
\caption{Center split: prediction performance for $S_{\mathrm{mom}}$ as a function of prefix length $T$, reported up to $T\le 3T_{\min}=36$. The dashed line indicates $T_{\min}=12$ (median minimum-location index of the good region), and the dotted line indicates $3T_{\min}$.}
\label{fig:center_metric_curves}
\end{figure}

\paragraph{Accuracy vs.\ horizon.}
Figure~\ref{fig:center_metric_curves} summarizes MAE, RMSE, and $R^2$ as functions of $T$ for all considered regressors.
A consistent trend is observed: prediction quality improves rapidly with increasing horizon and then gradually saturates.
Notably, strong predictive performance is achieved already around the characteristic minimum-location scale ($T\approx T_{\min}$), supporting the feasibility of early reliability estimation before the full profile minimum is necessarily observed.
Across horizons, nonlinear models (kNN and tree ensembles) substantially outperform linear baselines (Ridge/ElasticNet), especially for short prefixes, indicating that the mapping from early proxy dynamics to $S_{\mathrm{mom}}$ is strongly nonlinear.

\paragraph{Best model per horizon.}
Table~\ref{tab:center_best_by_T} lists the best-performing model (by $R^2$) for each scanned horizon $T\in\{1,3,5,\ldots,35\}$, together with the corresponding fit time and test-time inference costs.
The selection varies with horizon: Gradient Boosting is preferred at the very smallest horizons, Random Forest dominates at intermediate horizons, and kNN becomes competitive at larger horizons.
Overall, test-time costs remain small (microsecond to tens-of-microseconds per sample at this grid size), supporting the use of the learned diagnostic layer for early decision support.

\begin{table}[t]
\centering
\scriptsize
\setlength{\tabcolsep}{5pt}
\renewcommand{\arraystretch}{1.1}
\caption{Center split: best model by $R^2$ at each scanned horizon $T$. $R^2$ is rounded to three decimals; timing values are reported in seconds (scientific notation, two decimals).}
\label{tab:center_best_by_T}
\begin{tabular}{r l r r r r}
\hline
$T$ & best model & $R^2$ & fit time (s) & test time (s) & test per-sample (s) \\
\hline
1  & \texttt{grad\_boost}    & 0.476 & 2.91e-01 & 3.06e-03 & 2.12e-06 \\
3  & \texttt{grad\_boost}    & 0.670 & 6.30e-01 & 3.05e-03 & 2.12e-06 \\
5  & \texttt{random\_forest} & 0.791 & 3.12e-01 & 5.00e-02 & 3.47e-05 \\
7  & \texttt{random\_forest} & 0.836 & 3.46e-01 & 9.87e-02 & 6.85e-05 \\
9  & \texttt{random\_forest} & 0.873 & 3.67e-01 & 4.93e-02 & 3.42e-05 \\
11 & \texttt{random\_forest} & 0.893 & 4.48e-01 & 4.90e-02 & 3.40e-05 \\
13 & \texttt{random\_forest} & 0.909 & 4.48e-01 & 4.93e-02 & 3.43e-05 \\
15 & \texttt{knn}            & 0.928 & 5.41e-03 & 3.76e-02 & 2.61e-05 \\
17 & \texttt{knn}            & 0.934 & 8.12e-04 & 2.68e-02 & 1.86e-05 \\
19 & \texttt{knn}            & 0.942 & 1.04e-03 & 2.36e-02 & 1.64e-05 \\
21 & \texttt{knn}            & 0.942 & 1.18e-03 & 2.42e-02 & 1.68e-05 \\
23 & \texttt{knn}            & 0.944 & 1.08e-03 & 2.28e-02 & 1.59e-05 \\
25 & \texttt{knn}            & 0.948 & 1.10e-03 & 2.33e-02 & 1.62e-05 \\
27 & \texttt{knn}            & 0.946 & 1.37e-03 & 2.66e-02 & 1.85e-05 \\
29 & \texttt{grad\_boost}    & 0.946 & 2.87e-01 & 3.28e-03 & 2.28e-06 \\
31 & \texttt{knn}            & 0.955 & 1.17e-03 & 2.35e-02 & 1.63e-05 \\
33 & \texttt{knn}            & 0.956 & 1.19e-03 & 2.41e-02 & 1.67e-05 \\
35 & \texttt{knn}            & 0.956 & 1.18e-03 & 2.39e-02 & 1.66e-05 \\
\hline
\end{tabular}
\end{table}

\begin{figure*}[t]
\centering
\begin{minipage}{0.49\linewidth}
\centering
\includegraphics[width=\linewidth]{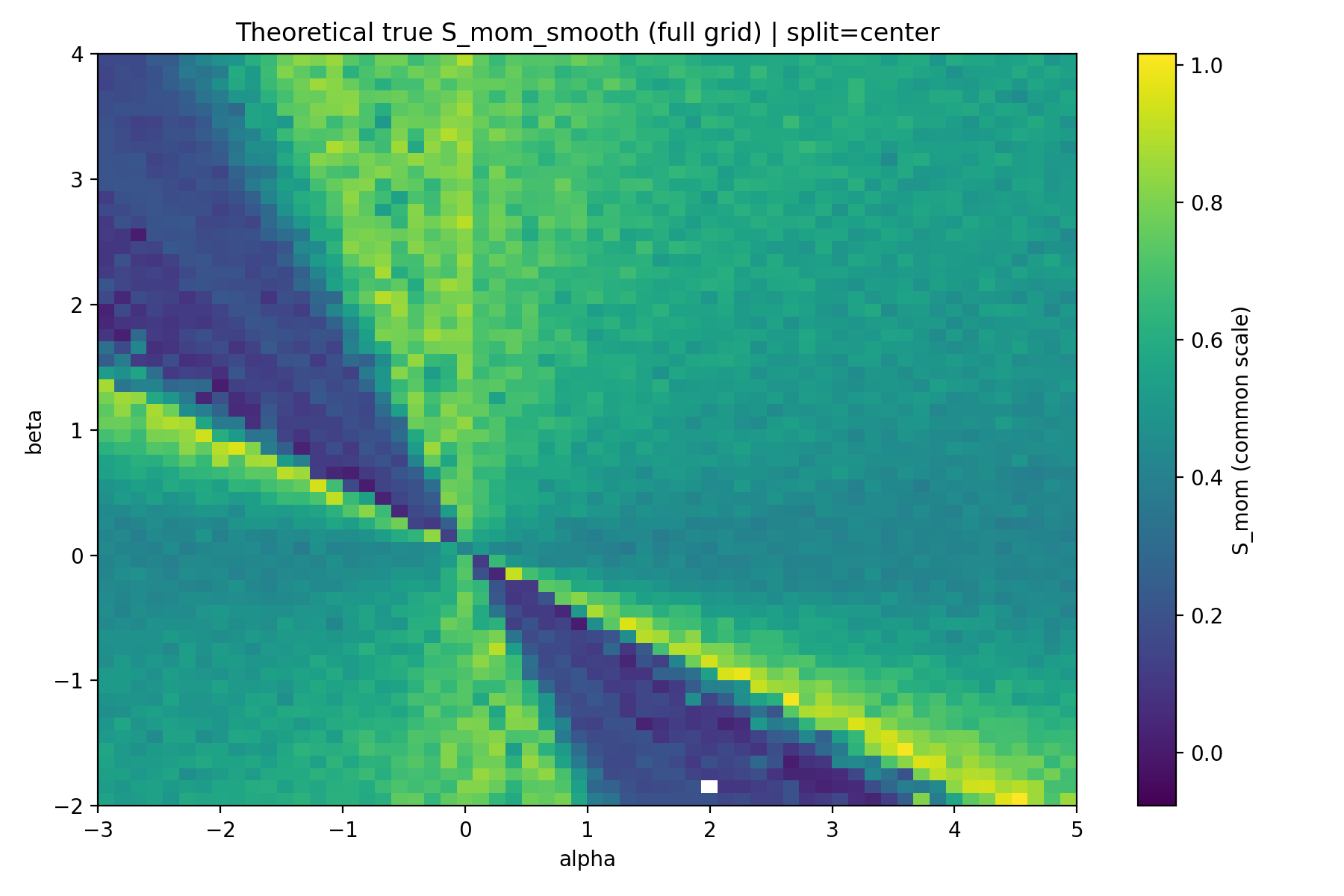}\\[-2pt]
\textbf{(a)} Theoretical $S_{\mathrm{mom}}$ (full grid).
\end{minipage}\hfill
\begin{minipage}{0.49\linewidth}
\centering
\includegraphics[width=\linewidth]{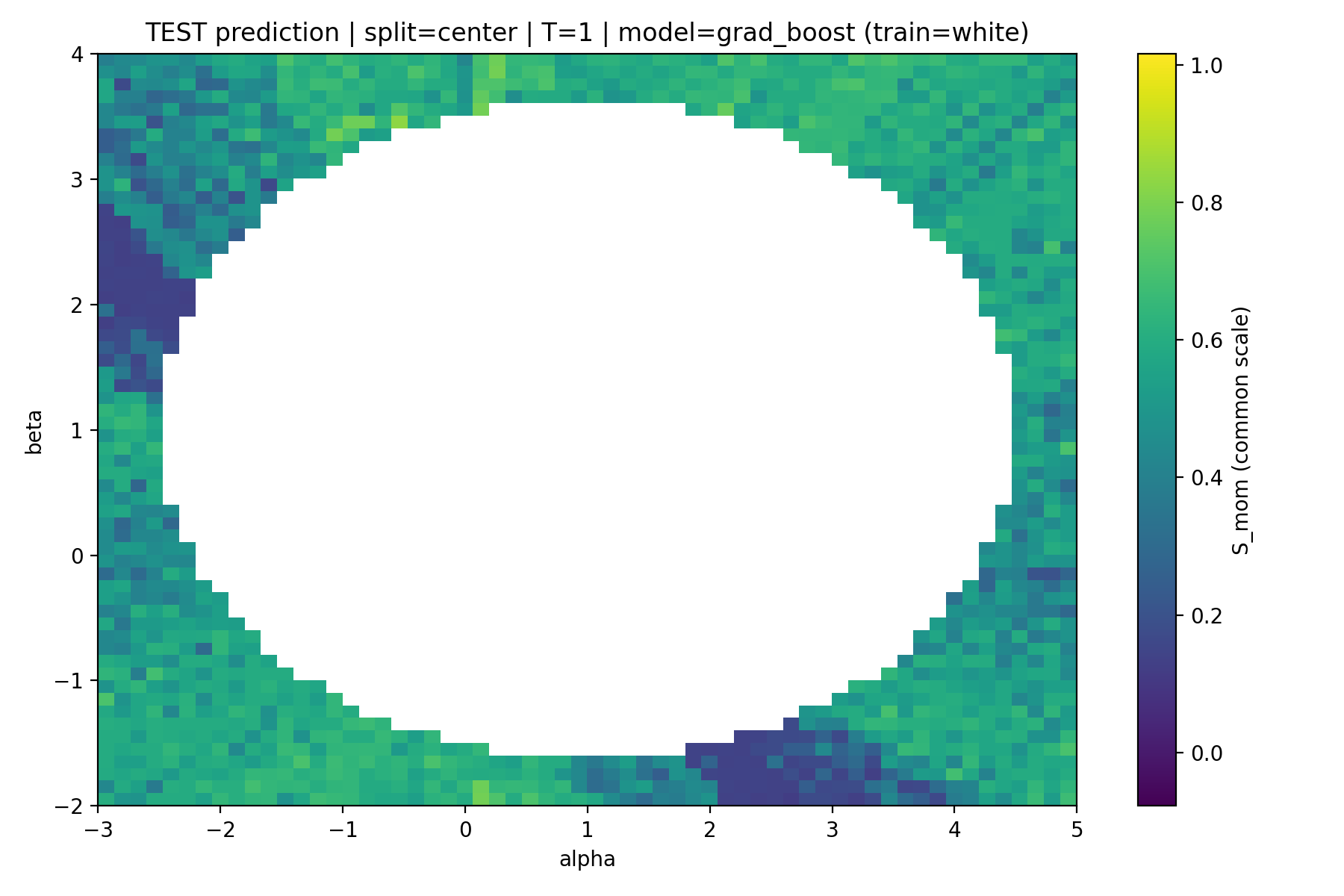}\\[-2pt]
\textbf{(b)} Test prediction, $T=1$ (\texttt{grad\_boost}).
\end{minipage}

\vspace{6pt}

\begin{minipage}{0.49\linewidth}
\centering
\includegraphics[width=\linewidth]{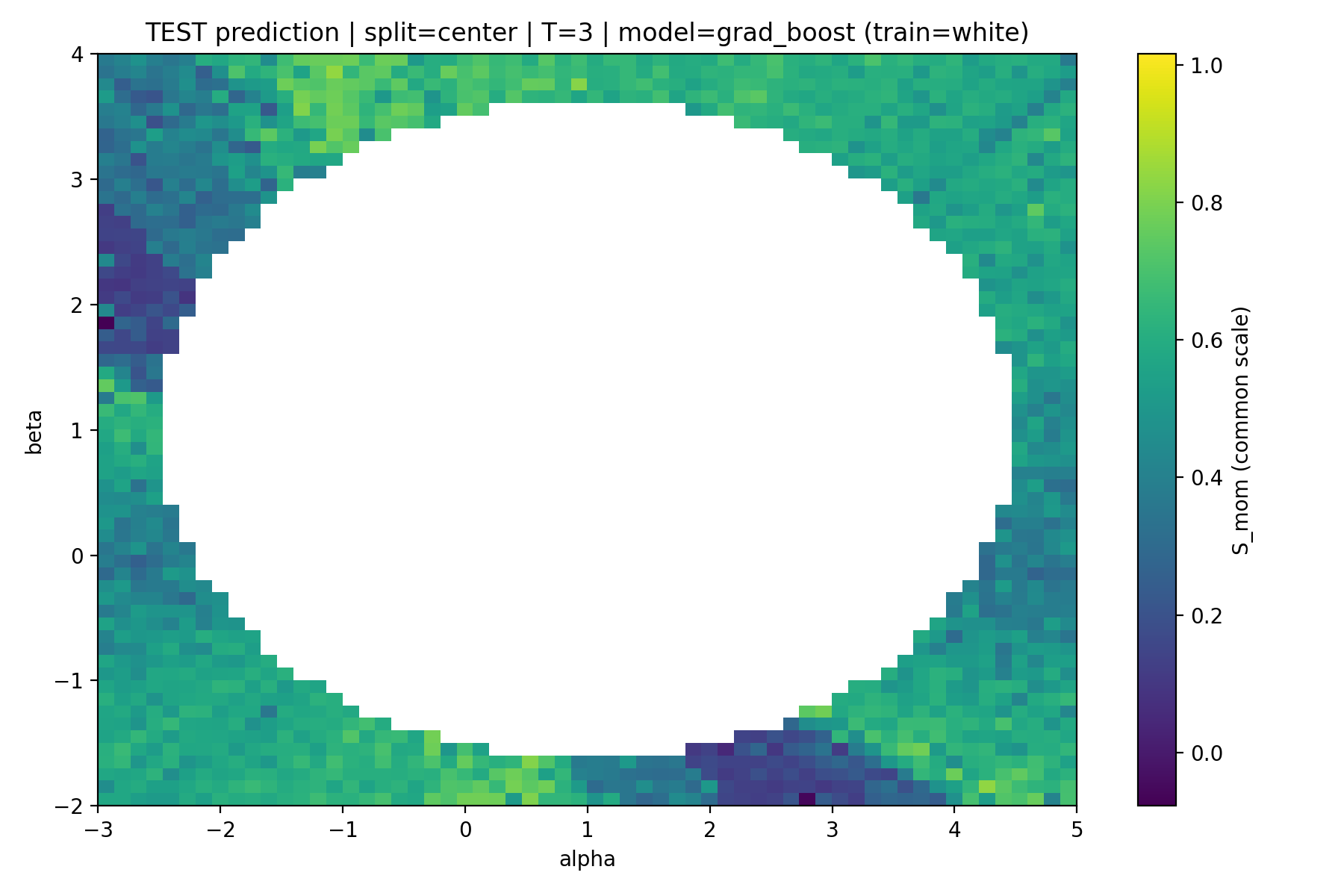}\\[-2pt]
\textbf{(c)} Test prediction, $T=3$ (\texttt{grad\_boost}).
\end{minipage}\hfill
\begin{minipage}{0.49\linewidth}
\centering
\includegraphics[width=\linewidth]{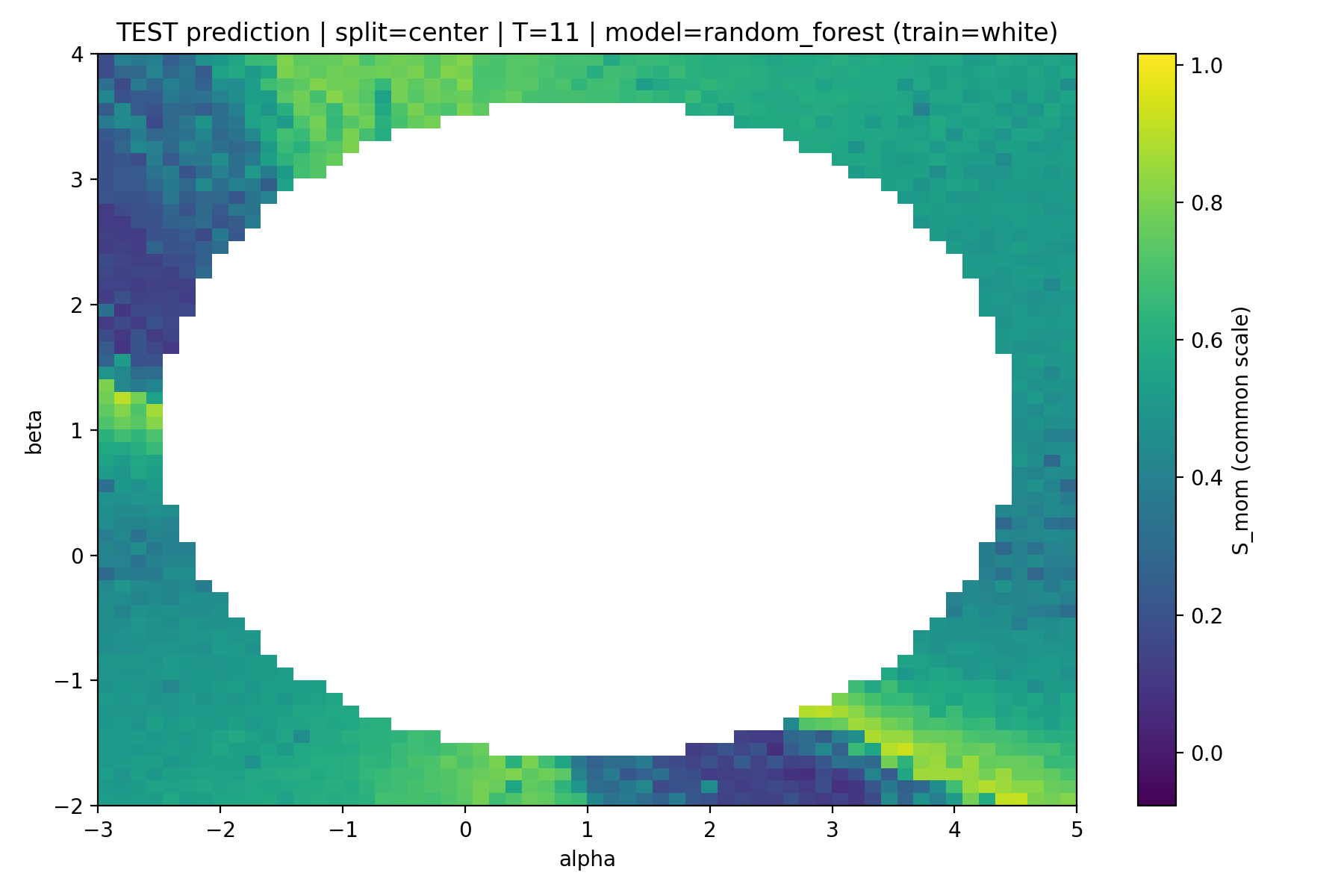}\\[-2pt]
\textbf{(d)} Test prediction, $T=11$ (\texttt{random\_forest}).
\end{minipage}

\caption{Center split: heatmap-level comparison between the theoretical $S_{\mathrm{mom}}$ landscape (a) and test-only predictions at increasing horizons (b)--(d). White region indicates training points (not plotted), while colored pixels correspond to test-region predictions on a common scale.}
\label{fig:center_heatmap_pred}
\end{figure*}

\paragraph{Heatmap-level comparison (test-only predictions).}
To visualize spatial generalization, Figure~\ref{fig:center_heatmap_pred} compares the theoretical $S_{\mathrm{mom}}$ landscape against test-only predicted heatmaps at representative horizons $T\in\{1,3,11\}$ using the best model selected at each horizon (Table~\ref{tab:center_best_by_T}). The training region is shown as white to emphasize out-of-sample behavior.
As $T$ increases, the predicted periphery increasingly recovers the main geometric structures of the true map, consistent with the monotone gains observed in Figure~\ref{fig:center_metric_curves}.

\subsection{Multi-Horizon Regression Results (Random Split)}
\begin{figure}[t]
\centering
\begin{minipage}{0.32\linewidth}
\centering
\includegraphics[width=\linewidth]{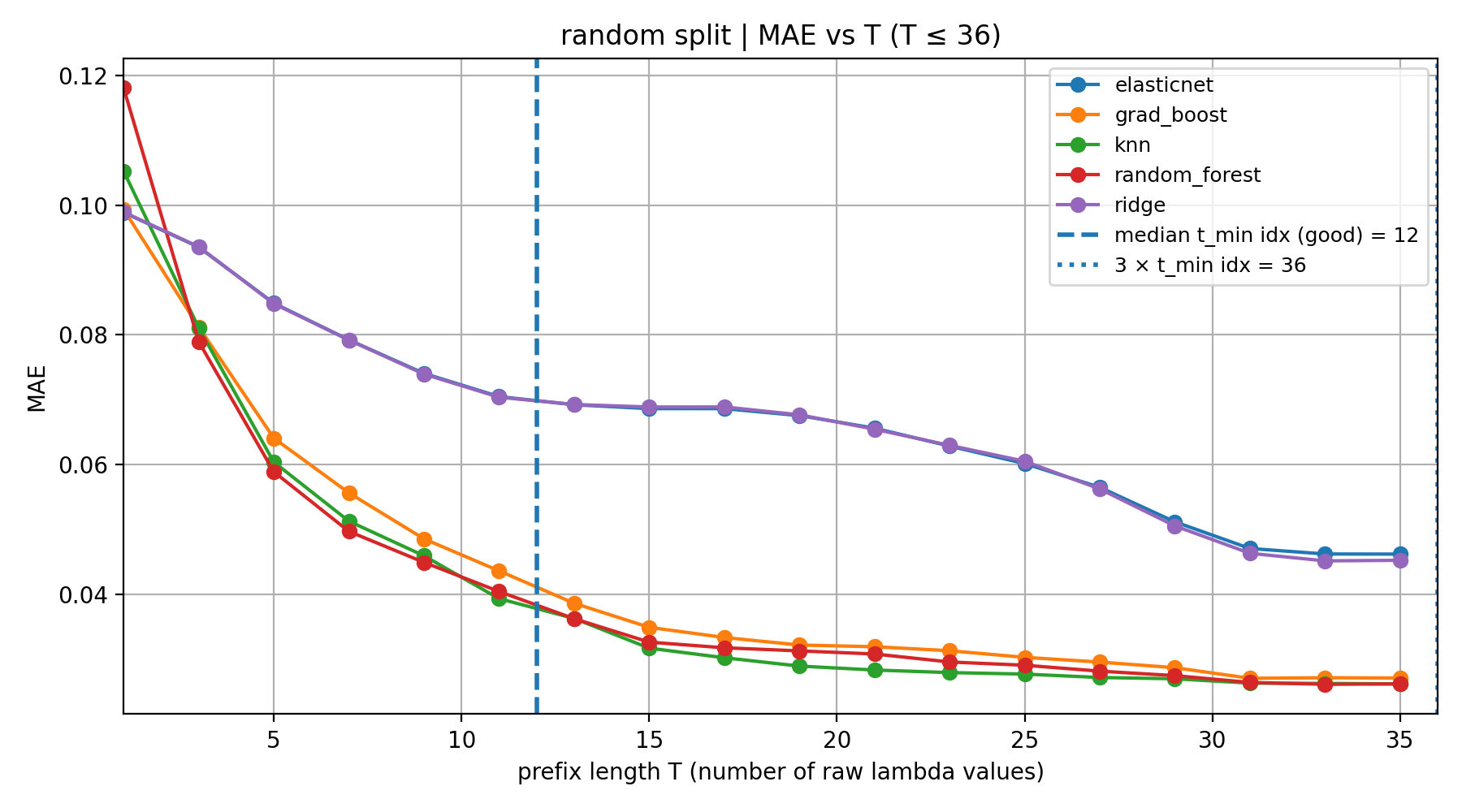}\\[-2pt]
(a) MAE vs.\ $T$.
\end{minipage}\hfill
\begin{minipage}{0.32\linewidth}
\centering
\includegraphics[width=\linewidth]{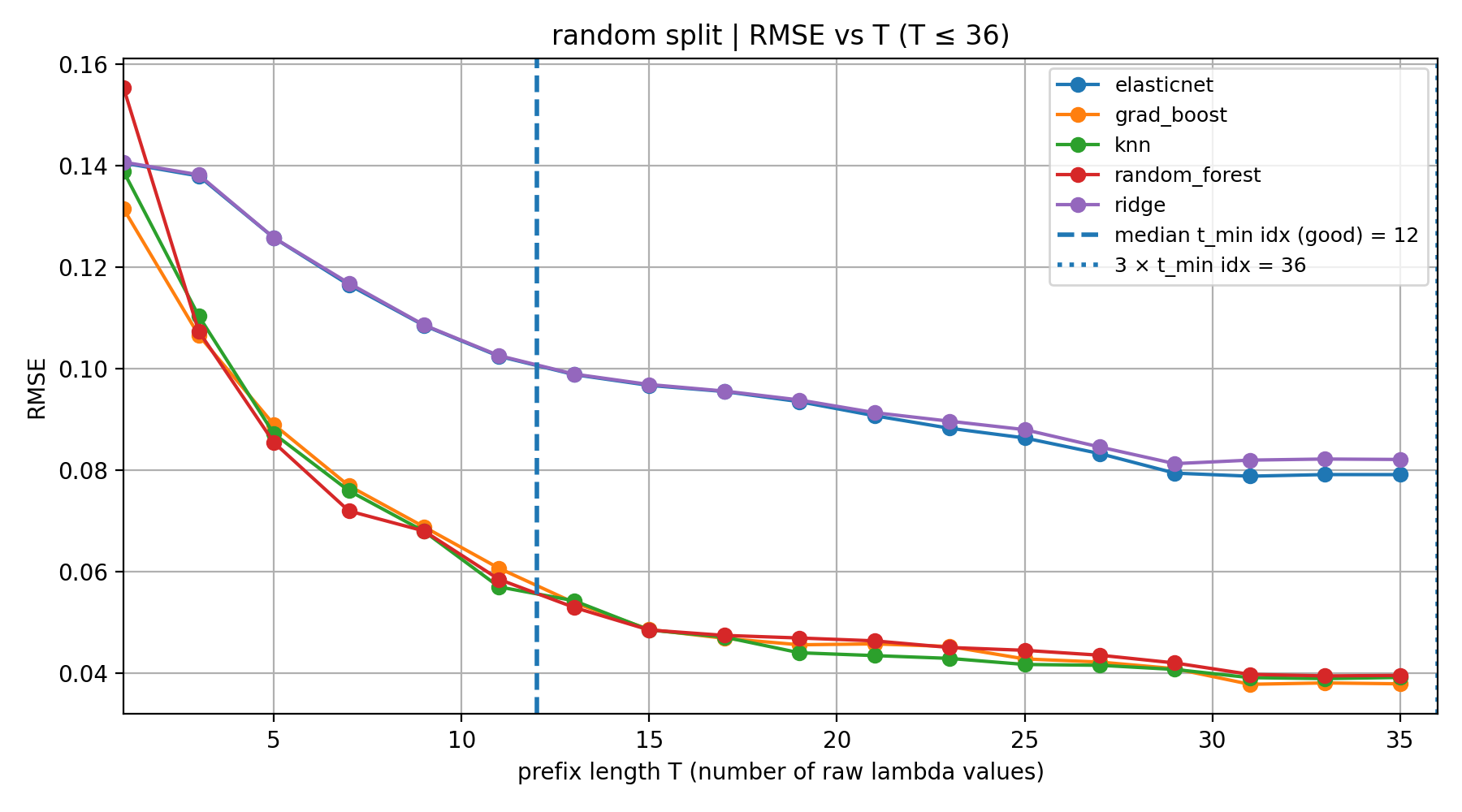}\\[-2pt]
(b) RMSE vs.\ $T$.
\end{minipage}\hfill
\begin{minipage}{0.32\linewidth}
\centering
\includegraphics[width=\linewidth]{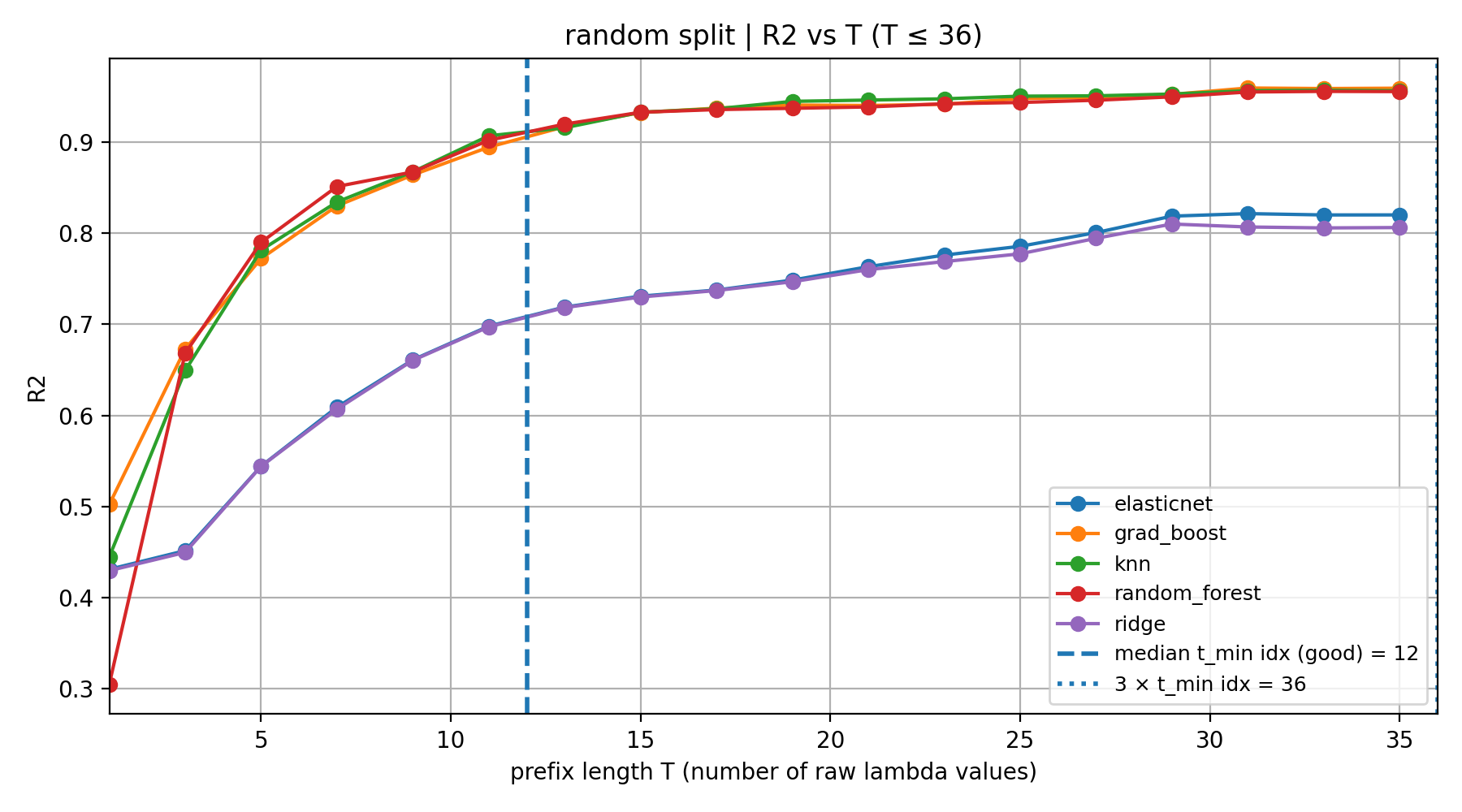}\\[-2pt]
(c) $R^2$ vs.\ $T$.
\end{minipage}
\caption{Random split: prediction performance for $S_{\mathrm{mom}}$ as a function of prefix length $T$, reported up to $T\le 3T_{\min}=36$. The dashed line indicates $T_{\min}=12$, and the dotted line indicates $3T_{\min}$.}
\label{fig:random_metric_curves}
\end{figure}

Under the random split, training and test samples are drawn from the same global distribution over the parameter grid. Figure~\ref{fig:random_metric_curves} reports MAE, RMSE, and $R^2$ as functions of horizon $T$, again truncated at $T\le 3T_{\min}=36$ for visual comparability with the center-split analysis.
The overall behavior mirrors that of the center-to-periphery split: prediction quality improves sharply at small horizons and then gradually saturates. The absolute metrics are slightly better than in the center split, as expected in a less extrapolative setting, but the qualitative pattern remains the same.

\paragraph{Best model per horizon.}
Table~\ref{tab:random_best_by_T} lists the best-by-$R^2$ model at each scanned horizon under the random split. The same main pattern appears as in the center split: Gradient Boosting is strongest at the very shortest horizons, while Random Forest and kNN dominate most larger horizons. Test-time inference remains negligible.

\begin{table}[t]
\centering
\scriptsize
\setlength{\tabcolsep}{4pt}
\renewcommand{\arraystretch}{1.05}
\caption{Random split: best model by $R^2$ at each scanned horizon $T$. $R^2$ is rounded to three decimals; timing values are reported in seconds (scientific notation, two decimals).}
\label{tab:random_best_by_T}
\begin{tabular}{r l r r r r}
\hline
$T$ & best model & $R^2$ & fit (s) & test (s) & test/sample (s) \\
\hline
1  & \texttt{grad\_boost}    & 0.486 & 3.00e-01 & 3.41e-03 & 2.37e-06 \\
3  & \texttt{grad\_boost}    & 0.673 & 4.57e-01 & 3.44e-03 & 2.39e-06 \\
5  & \texttt{random\_forest} & 0.791 & 3.00e-01 & 5.00e-02 & 3.47e-05 \\
7  & \texttt{knn}            & 0.852 & 2.50e-03 & 2.61e-03 & 1.81e-06 \\
9  & \texttt{knn}            & 0.879 & 2.54e-03 & 2.54e-03 & 1.76e-06 \\
11 & \texttt{knn}            & 0.905 & 2.53e-03 & 2.50e-03 & 1.74e-06 \\
13 & \texttt{random\_forest} & 0.919 & 4.36e-01 & 1.01e-01 & 7.01e-05 \\
15 & \texttt{knn}            & 0.931 & 2.44e-03 & 2.42e-03 & 1.68e-06 \\
17 & \texttt{grad\_boost}    & 0.933 & 9.75e-01 & 3.28e-03 & 2.28e-06 \\
19 & \texttt{knn}            & 0.935 & 2.45e-03 & 2.39e-03 & 1.66e-06 \\
21 & \texttt{knn}            & 0.939 & 2.47e-03 & 2.44e-03 & 1.69e-06 \\
23 & \texttt{knn}            & 0.941 & 2.54e-03 & 2.62e-03 & 1.82e-06 \\
25 & \texttt{knn}            & 0.945 & 2.56e-03 & 2.44e-03 & 1.69e-06 \\
27 & \texttt{knn}            & 0.944 & 2.63e-03 & 2.56e-03 & 1.78e-06 \\
29 & \texttt{knn}            & 0.944 & 2.64e-03 & 2.45e-03 & 1.70e-06 \\
31 & \texttt{grad\_boost}    & 0.956 & 1.12e00 & 3.20e-03 & 2.22e-06 \\
33 & \texttt{grad\_boost}    & 0.957 & 1.13e00 & 3.25e-03 & 2.25e-06 \\
35 & \texttt{grad\_boost}    & 0.958 & 1.17e00 & 3.22e-03 & 2.24e-06 \\
\hline
\end{tabular}
\end{table}

\begin{figure*}[t]
\centering
\begin{minipage}{0.49\linewidth}
\centering
\includegraphics[width=\linewidth]{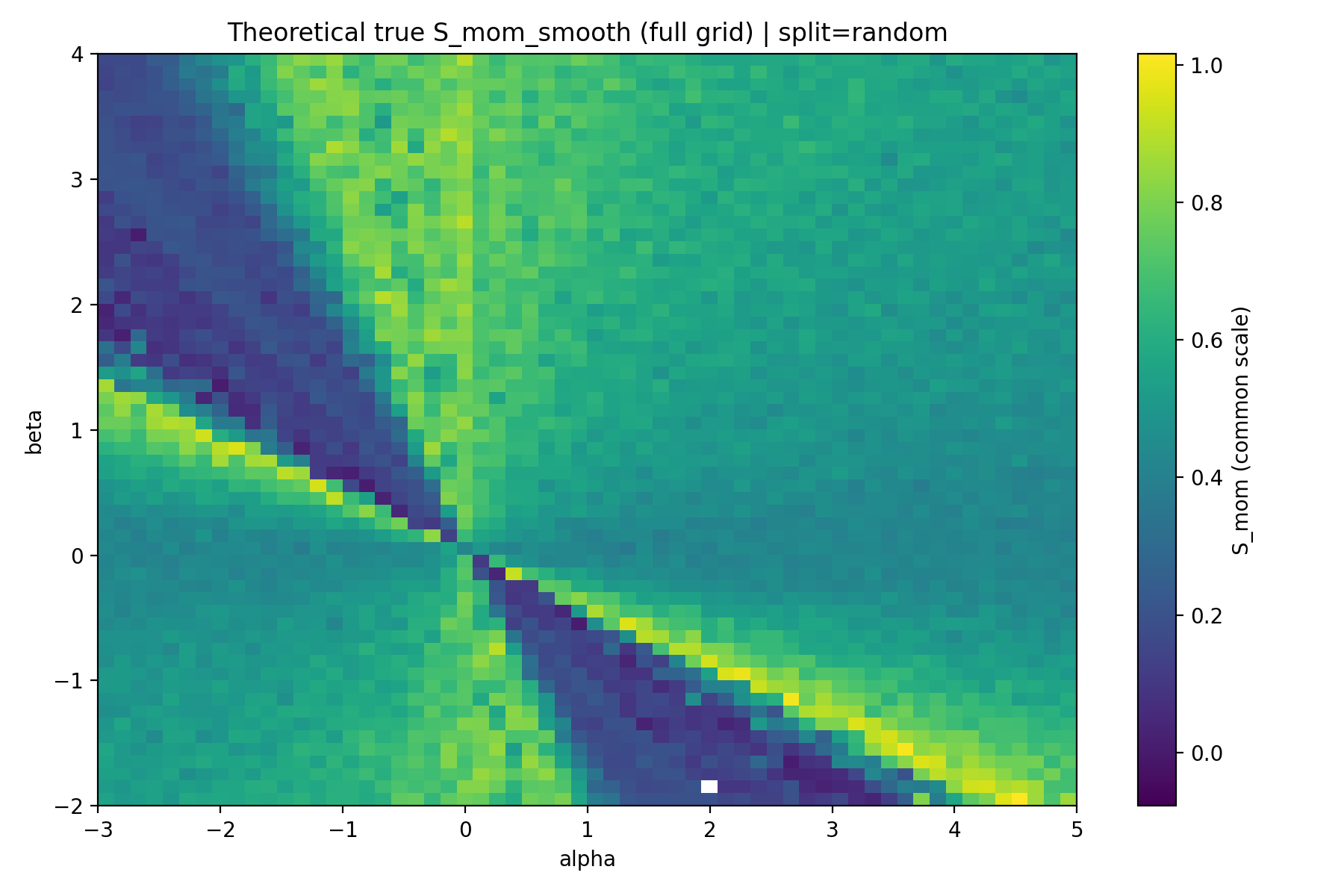}\\[-2pt]
\textbf{(a)} Theoretical $S_{\mathrm{mom}}$ (full grid).
\end{minipage}\hfill
\begin{minipage}{0.49\linewidth}
\centering
\includegraphics[width=\linewidth]{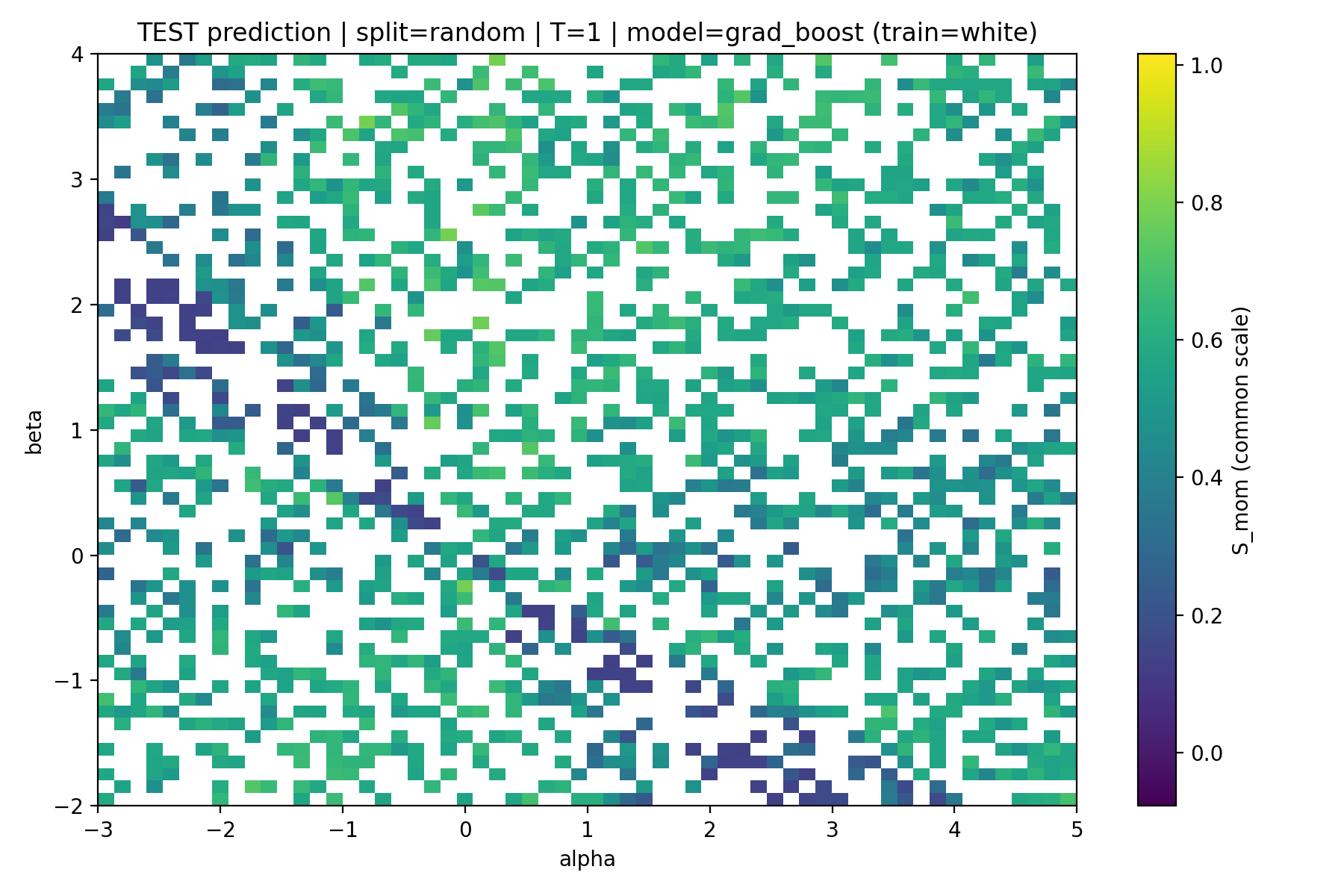}\\[-2pt]
\textbf{(b)} Test-only prediction, $T=1$ (\texttt{grad\_boost}; train=white).
\end{minipage}

\vspace{6pt}

\begin{minipage}{0.49\linewidth}
\centering
\includegraphics[width=\linewidth]{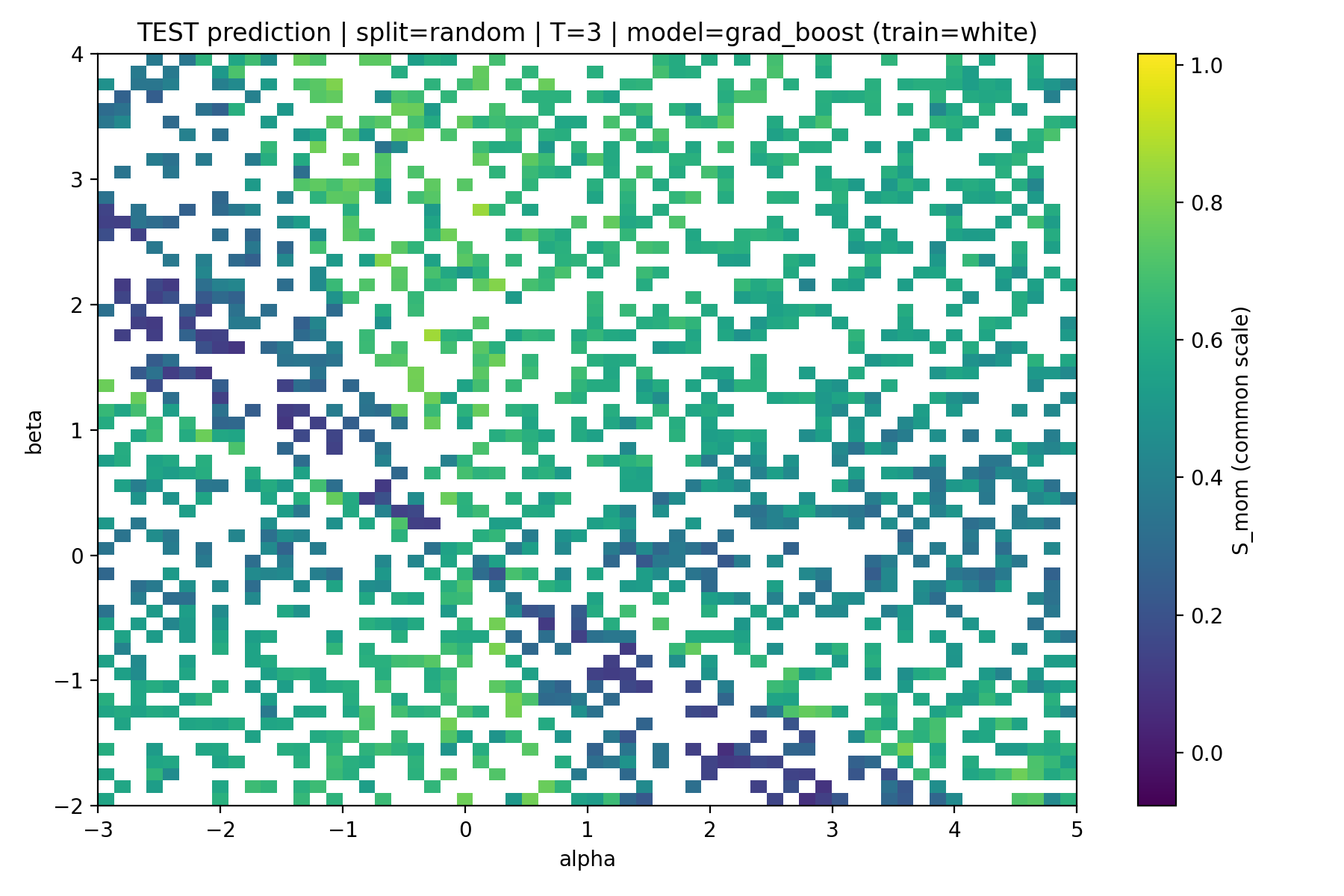}\\[-2pt]
\textbf{(c)} Test-only prediction, $T=3$ (\texttt{grad\_boost}; train=white).
\end{minipage}\hfill
\begin{minipage}{0.49\linewidth}
\centering
\includegraphics[width=\linewidth]{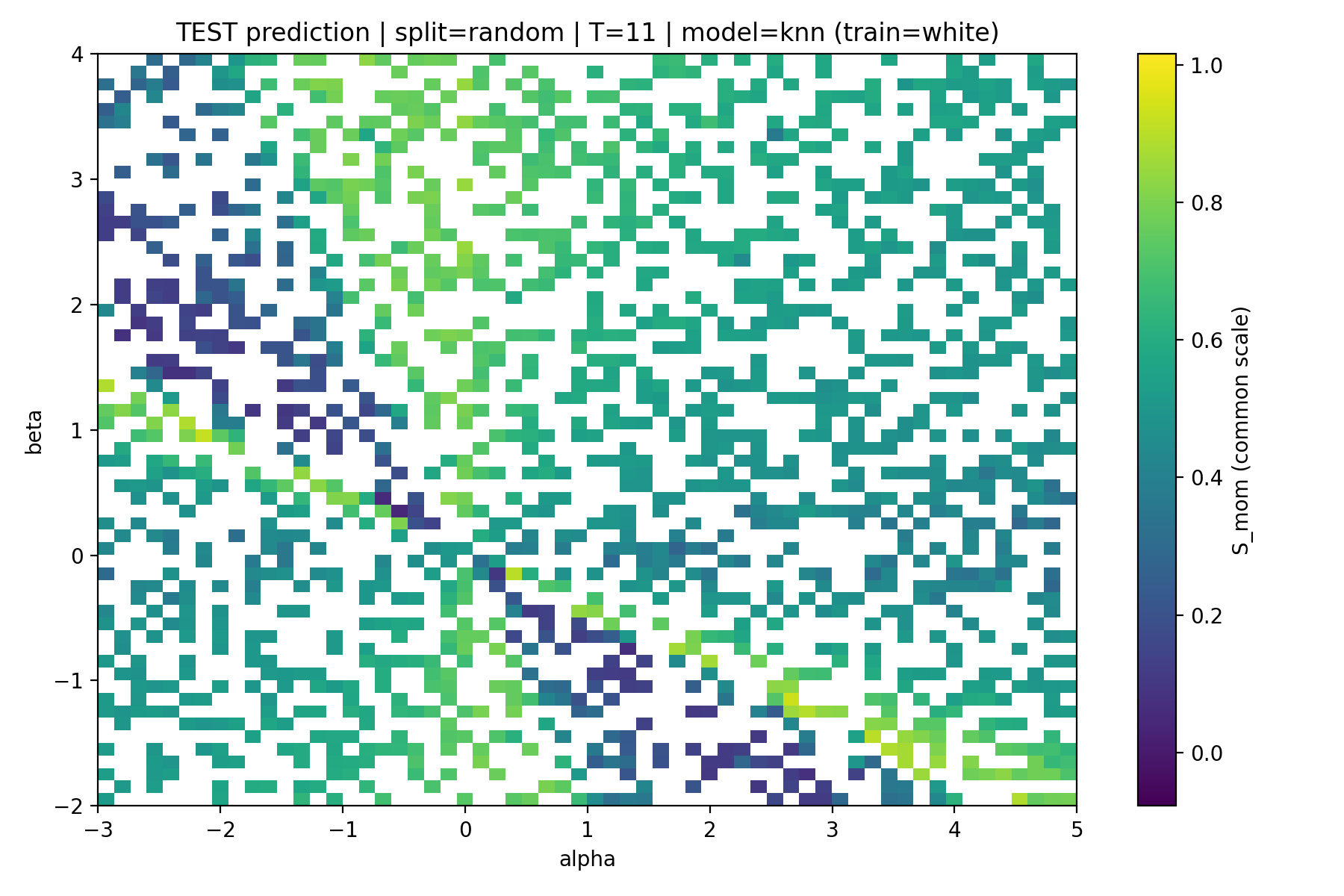}\\[-2pt]
\textbf{(d)} Test-only prediction, $T=11$ (\texttt{knn}; train=white).
\end{minipage}

\caption{Random split: heatmap-level comparison between the theoretical $S_{\mathrm{mom}}$ landscape (a) and test-only predictions at increasing horizons (b)--(d). White cells indicate training points (not plotted), while colored cells correspond to test-region predictions on a common scale.}
\label{fig:random_heatmap_pred}
\end{figure*}

\paragraph{Heatmap-level comparison (test-only predictions).}
Figure~\ref{fig:random_heatmap_pred} presents test-only predicted heatmaps for representative horizons $T\in\{1,3,11\}$ under the random split, using the corresponding best models from Table~\ref{tab:random_best_by_T}. Because the train/test partition is spatially interspersed, the training region appears as a punctured white mask rather than a contiguous block. Nevertheless, as $T$ increases, the predicted values on the held-out points increasingly reproduce the true landscape of $S_{\mathrm{mom}}$.

\subsection{Computational Benefit of Early Prediction}
The practical role of the predictive layer is not to compute the roots themselves, nor to optimize the initial guesses, but to screen parameter configurations $(\alpha,\beta)$ by predicting the final reliability score $S_{\mathrm{mom}}$ from an early prefix of the proxy profile. This distinction is important: the computational saving arises from avoiding the cost of generating the full diagnostic profile for every parameter pair.

In the present setup, the full ground-truth diagnostic at a single grid point uses $K=200$ solver iterations per trajectory. By contrast, a predictor operating at prefix length $T$ requires only the initial portion of the proxy profile. Because each proxy value $\lambda_1(t_{\mathrm{end}})$ is itself built from a delay embedding of length $L=5$ and forecast horizons up to $h_{\max}=5$, the horizon index $T$ does not correspond one-to-one to the raw solver iteration count. A simple accounting gives an approximate requirement of
\begin{equation}
k_{\mathrm{req}}(T) \approx T + L + h_{\max} - 1 = T + 9,
\end{equation}
solver iterations per trajectory in the current configuration. Therefore, the relative diagnostic cost at horizon $T$ is approximately $(T+9)/200$, and the corresponding speedup with respect to the full diagnostic is about $200/(T+9)$.

Table~\ref{tab:cost_vs_quality} combines this cost estimate with the observed predictive quality at representative horizons. The results show that useful accuracy is achieved very early. In particular, at $T=11$, which is close to the characteristic minimum-location scale identified from the good-region timing statistics, the predictive layer already achieves $R^2\approx0.89$--$0.91$ while using only about $20$ solver iterations instead of the full $200$. This corresponds to an approximately tenfold reduction in diagnostic cost per trajectory. Even the coarser horizon $T=3$ requires only about $12$ solver iterations, yielding an approximately $16.7\times$ reduction in solver budget while still recovering a useful fraction of the final reliability landscape.

\begin{table}[t]
\centering
\caption{Approximate diagnostic cost versus predictive quality at representative horizons. The solver-iteration requirement is estimated as $k_{\mathrm{req}}(T)\approx T+9$ for the present proxy construction ($L=5$, $h_{\max}=5$), and the full diagnostic uses $K=200$ iterations per trajectory.}
\label{tab:cost_vs_quality}
\begin{tabular}{c c c c c}
\hline
$T$ & approx. solver iterations & speedup vs. full $K=200$ & center split $R^2$ & random split $R^2$ \\
\hline
1  & 10 & $20.0\times$ & 0.476 & 0.486 \\
3  & 12 & $16.7\times$ & 0.670 & 0.673 \\
11 & 20 & $10.0\times$ & 0.893 & 0.905 \\
35 & 44 & $4.5\times$ & 0.956 & 0.958 \\
\hline
\end{tabular}
\end{table}

This cost reduction is particularly relevant for screening unfavorable parameter configurations. The timing histograms in Figure~\ref{fig:timing_hists} show that the characteristic profile events used by the full diagnostic, such as the minimum of $\tilde{\lambda}_1(t)$ and the first entry into the negative regime, are not uniformly early across the parameter grid. In less favorable regions, these events can occur substantially later than in the good subset, so waiting for the complete profile delays the rejection of poor parameter choices. By contrast, the learned predictor can issue an informative estimate of the eventual $S_{\mathrm{mom}}$ score from a much shorter prefix. In this sense, the proposed AI layer does not accelerate the underlying root-finding iterations themselves; rather, it accelerates the reliability assessment and parameter screening process built on top of the solver.

\section{Discussion}
\label{sec:discussion}

The results demonstrate that short prefixes of kNN--LLE proxy profiles contain substantial information about the final reliability score $S_{\mathrm{mom}}$, making early prediction feasible across a structured two-parameter solver landscape. Beyond the raw regression metrics, several points deserve further interpretation.

\paragraph{Why early prediction works.}
The proxy profile is designed to summarize local contractivity and expansivity trends in the iteration dynamics. The empirical curves in Figures~\ref{fig:center_metric_curves} and \ref{fig:random_metric_curves} exhibit a common pattern: a steep accuracy gain at small prefix lengths $T$, followed by a gradual saturation. This behavior is consistent with the role of $S_{\mathrm{mom}}$ as an early-and-strong contractivity indicator: once the initial portion of the proxy dynamics reveals the onset (or absence) of contractive tendencies, additional late-profile information yields diminishing improvements.

\paragraph{Interpretability of the predictive layer.}
An important feature of the proposed framework is that the machine learning models do not operate directly on raw solver iterates or full iteration trajectories. Instead, they use early prefixes of the kNN--LLE proxy profile, which encode interpretable local stability signatures of the underlying dynamics. In this sense, the predictive layer remains connected to the dynamical interpretation of solver behavior: the learned mapping is not from arbitrary high-dimensional traces to a score, but from contractive/expansive profile patterns to the final reliability indicator $S_{\mathrm{mom}}$. This makes the approach more transparent than a purely black-box surrogate built directly on raw iteration data.

\paragraph{Prediction before the characteristic minimum scale ($T<T_{\min}$).}
A key practical regime is prediction before the typical minimum-location index of the good region, $T_{\min}=12$. This corresponds to forecasting reliability before the visually apparent minimum of the smoothed proxy profile is necessarily observed. In both splits, strong accuracy is already achieved around $T\approx T_{\min}$, indicating that reliable discrimination is possible from early signatures alone. This is particularly important for early decision support, where the goal is to avoid spending the full computational budget required to reach the minimum if an unfavorable configuration can be detected earlier.

\paragraph{The extreme early-exit case $T=1$.}
The $T=1$ setting is especially notable: the model receives only the first raw proxy value $\lambda_1^{\mathrm{raw}}(1)$. Despite being a single scalar, this quantity is not a trivial instantaneous measurement. In the underlying kNN--LLE construction, each $\lambda_1(t_{\mathrm{end}})$ value is produced from a delay-embedding of length $L$ and a multi-horizon forecast over $h\in\{1,\ldots,5\}$, so the earliest proxy estimate already aggregates information from a short history and multiple forward steps. In that sense, $T=1$ still represents a compressed summary of early solver dynamics rather than a single iteration of the base root-finding scheme. Nevertheless, the errors at $T=1$ remain larger than for $T\ge 3$, and the heatmap comparisons confirm that very short horizons capture only coarse structures. This suggests that while a one-step diagnostic is feasible, a modest prefix length ($T\approx 3$--$5$) provides a substantially more reliable and stable prediction regime.

\paragraph{Random vs.\ center-to-periphery generalization.}
The two splitting strategies yield broadly comparable quantitative accuracy, but they probe different notions of generalization. The random split evaluates within-distribution performance, and the associated test-only heatmaps appear spatially punctured because the training region is interspersed throughout the domain. In contrast, the center-to-periphery split enforces a structured extrapolation test: the models are trained on central parameter regimes and evaluated on boundary regimes, which is more representative of out-of-distribution deployment scenarios (e.g., screening untested parameter regions or transferring to nearby regimes). The fact that the center split still achieves high $R^2$ at moderate horizons supports the robustness of the learned mapping across the parameter domain, although the visual recovery of fine geometric details is naturally slower than in the i.i.d.\ case.

\paragraph{Model behavior and inductive biases.}
Across both splits, nonlinear regressors (kNN, Gradient Boosting, Random Forest) consistently dominate the linear baselines (Ridge and Elastic Net), particularly at short horizons. This indicates that the mapping from early proxy dynamics to $S_{\mathrm{mom}}$ is strongly nonlinear and regime-dependent. Ridge and Elastic Net provide stable baselines under correlated features, but their limited expressivity is reflected by systematically lower $R^2$ and their absence from best-by-$T$ selections. Among nonlinear models, the differences are comparatively small on the metric curves, suggesting that several model families can approximate the target mapping well once sufficient prefix information is available. The best-by-$T$ tables show that model preference can shift with horizon: boosting methods can be advantageous at extremely small horizons, while kNN or Random Forest becomes competitive as $T$ increases. Importantly, regardless of the chosen model, inference costs remain negligible relative to the profile-generation pipeline, reinforcing the viability of using the predictor as a lightweight reliability layer.

\paragraph{Heatmap-level interpretation.}
The heatmap comparisons provide an interpretable spatial validation beyond scalar metrics. For small horizons, predicted maps reproduce only large-scale features (e.g., broad high-score ridges and low-score wedges), while increasing $T$ progressively restores the geometry of the theoretical landscape. This qualitative improvement is consistent with the monotone gains in MAE/RMSE and $R^2$. Under center-to-periphery splitting, the predicted boundary region becomes visually close to the theoretical map at moderate horizons, which is the desired behavior for screening unseen parameter regions.

\paragraph{Generality beyond a specific equation and open questions.}
A central motivation of the proposed diagnostic layer is that it relies on observable dynamics (proxy contractivity profiles) rather than closed-form properties of a particular nonlinear equation. Consequently, the approach is expected to be less equation-specific than direct surrogate modeling of solver outputs. In principle, a model trained on one problem class may transfer to another, provided the early proxy profiles encode comparable stability signatures; however, this transferability requires dedicated validation. Future work should therefore test cross-problem generalization (training on one equation/model family and evaluating on another) and quantify how changes in profiling hyperparameters ($L$, horizon set, smoothing width) affect stability of the learned mapping and the best-by-$T$ selection.

\paragraph{Possible applications in systems biology and related steady-state studies.}
A particularly promising application perspective of the proposed framework arises in parameterized steady-state problems from systems biology and related areas. In biochemical reaction networks, metabolic ODE models, and signaling systems, biologically relevant steady states are often computed repeatedly across many parameter configurations for parameter estimation, sensitivity analysis, bifurcation studies, or regime screening \cite{lines2019steady,lakrisenko2023adjoint,rosenblatt2016constraints,fiedler2016steady}. In such settings, repeated nonlinear solves can become a major computational bottleneck, and convergence may depend strongly on the parameter regime and on the numerical behavior of the solver.

From this viewpoint, the present diagnostic layer could serve as an early screening tool for repeated steady-state computations. Rather than replacing the underlying nonlinear solver, it could help identify parameter regions that are numerically favorable, likely to yield stable steady states, or conversely prone to instability or wasted computational effort. This perspective is also consistent with recent literature on stability-constrained parameter exploration and parameter-space partitioning in biological models \cite{kariya2022stability,hernandez2023analytic,bradford2020multistationarity}. At the same time, the current paper does not yet demonstrate such an application directly; this should therefore be regarded as a realistic and literature-supported direction for future work rather than as a validated biological case study.

\paragraph{Summary of main observations.}
\begin{itemize}
\item Predictive accuracy improves rapidly with horizon $T$ and saturates thereafter; strong performance emerges near the characteristic scale $T_{\min}=12$, enabling early reliability prediction before the typical minimum is observed.
\item The $T=1$ regime is feasible but relatively noisy; modest horizons ($T\approx 3$--$5$) provide a more reliable early-exit operating point.
\item Nonlinear models (kNN and tree ensembles) substantially outperform linear baselines (Ridge/Elastic Net), indicating a strongly nonlinear mapping from early proxy dynamics to $S_{\mathrm{mom}}$.
\item Random and center-to-periphery splits yield broadly comparable accuracy, while the center split provides a more stringent test of out-of-distribution generalization across the parameter domain.
\item Heatmap-level comparisons confirm that increasing $T$ yields progressively more faithful recovery of the true spatial reliability landscape, supporting the interpretability and practical utility of the proposed diagnostic layer.
\end{itemize}

\section{Conclusion}
\label{sec:conclusion}
We presented an interpretable AI-assisted reliability diagnostic for a two-parameter parallel root-finding scheme, combining kNN--LLE proxy profiling with multi-horizon prediction of the profile-derived score $S_{\mathrm{mom}}$. The study shows that short early prefixes of the raw proxy profile already contain sufficient information to predict the final reliability score with high accuracy across the $(\alpha,\beta)$ parameter domain, including under the more demanding center-to-periphery split.

The main practical consequence is that reliability screening can be performed well before the full proxy profile is available. In the present experiments, prediction at $T=11$ achieves $R^2\approx0.89$--$0.91$ while requiring only about one tenth of the full solver budget needed to compute the complete diagnostic profile. This makes the proposed layer suitable for early rejection of unfavorable parameter configurations and for integrating lightweight decision support into solver workflows. Future work will focus on cross-problem transfer, robustness with respect to profiling hyperparameters, and adaptive policies that couple the learned diagnostic directly to runtime solver control.


\appendix
\section{Auxiliary Solver Validation of the Base Scheme}
\label{app:solver_validation}
This appendix collects solver-level validation results that support the numerical quality of the underlying parallel scheme but are not part of the main predictive contribution of the paper. All experiments in this appendix are performed at the fixed representative parameter pair $(\alpha,\beta)=(-0.1,4.0)$, chosen from the stable region of the $(\alpha,\beta)$ landscape.

To examine the robustness of the iterative scheme, three initialization strategies are used:

\textbf{Near–root:}
\[
z_k^{(0)}=\alpha_k+10^{-2}\varepsilon_k, \qquad \varepsilon_k\in[-1,1], \quad k=1,\ldots,7.
\]

\textbf{Moderate:}
\[
z_k^{(0)} = r\,e^{\frac{2\pi i (k-1)}{7}}, \qquad r=5,\; k=1,\ldots,7.
\]

\textbf{Random:}
\[
z_k^{(0)} \sim \mathcal{U}([-10,10] + i[-10,10]), \qquad k=1,\ldots,7.
\]

\subsection{Sensitivity to Initial Guesses}
\begin{figure}[!htbp]
\centering
\begin{subfigure}[b]{0.32\linewidth}
    \centering
    \includegraphics[width=\linewidth]{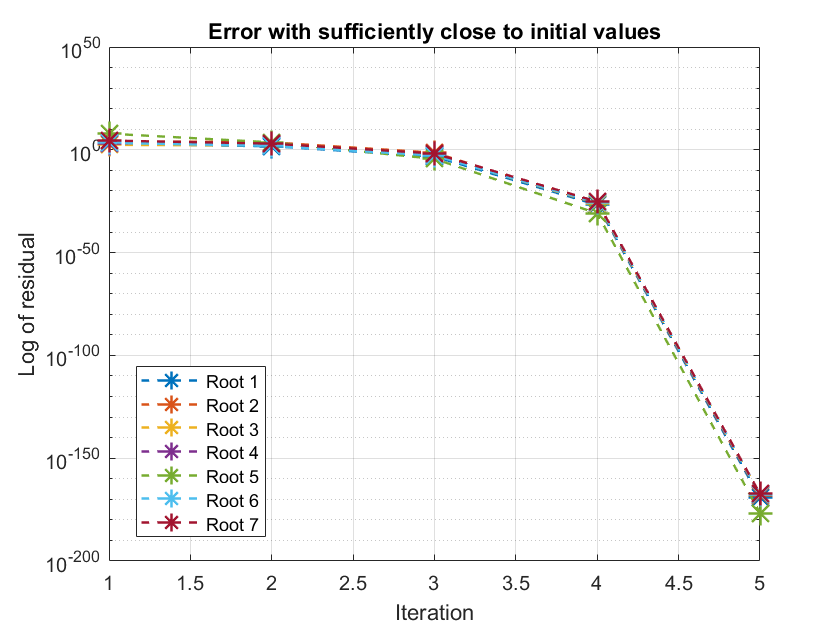}
    \caption{Near-root initialization}
    \label{fig:E1a}
\end{subfigure}\hfill
\begin{subfigure}[b]{0.32\linewidth}
    \centering
    \includegraphics[width=\linewidth]{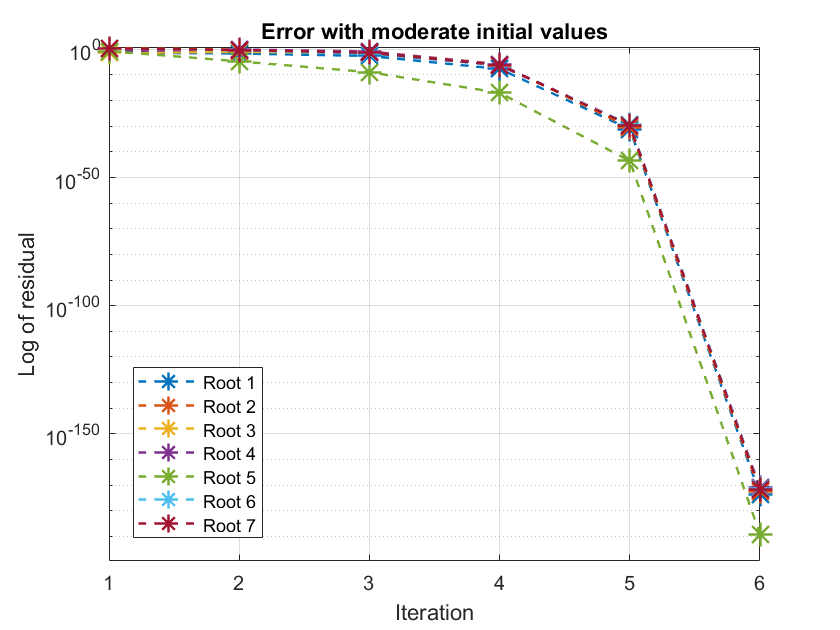}
    \caption{Moderate initialization}
    \label{fig:E1b}
\end{subfigure}\hfill
\begin{subfigure}[b]{0.32\linewidth}
    \centering
    \includegraphics[width=\linewidth]{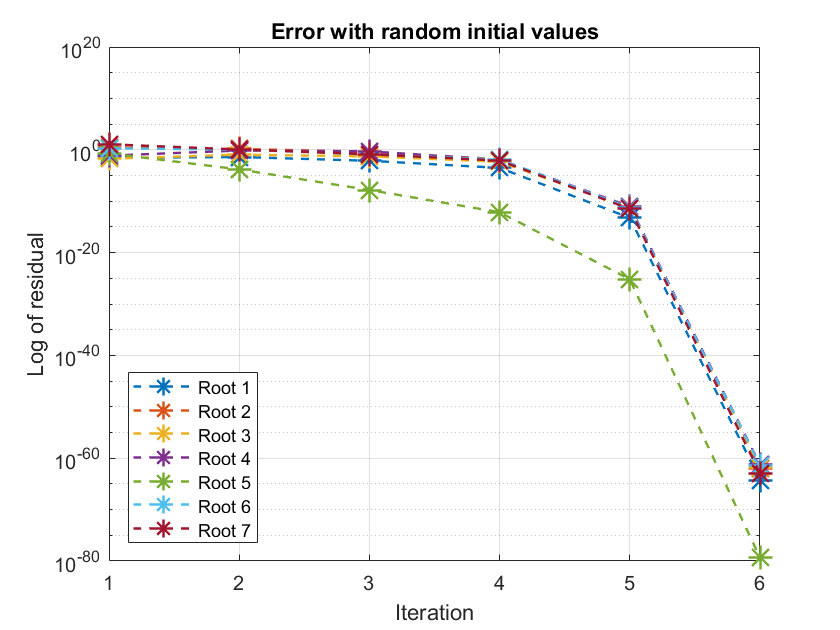}
    \caption{Random initialization}
    \label{fig:E1c}
\end{subfigure}
\caption{Early-iteration convergence of the proposed parallel scheme under different initialization scenarios, showing $\log_{10}(E^{(k)})$ versus iteration index.}
\label{E1}
\end{figure}

The convergence behavior of simultaneous root-finding methods depends on the quality of the initial approximations. To illustrate the robustness of the base solver under different starts, we consider three initialization scenarios: near-root initialization, moderately distant initialization, and random initialization. For each case, we monitor the maximum per-root error at iteration $k$,
\begin{equation}
E^{(k)}=\max_{1\le i\le n}\left|z_i^{(k)}-\xi_i\right|.
\end{equation}
The corresponding convergence histories are shown in Figure~\ref{E1}, while Table~\ref{tab:early_errors} summarizes the early-iteration error decay.

\begin{table}[!htbp]
\centering
\caption{Early iteration error behavior of the proposed scheme under different initialization scenarios.}
\begin{tabular}{c c c c c}
\hline
Iteration $k$ & Near-root init. & Moderate init. & Random init. & Order (approx.) \\
\hline
1 & $1.26\times10^{-2}$ & $3.35\times10^{-1}$ & $8.44\times10^{-1}$ & -- \\
2 & $2.86\times10^{-7}$ & $5.71\times10^{-5}$ & $3.47\times10^{-4}$ & 3.92 \\
3 & $4.66\times10^{-21}$ & $7.44\times10^{-9}$ & $9.78\times10^{-10}$ & 4.15 \\
4 & $2.16\times10^{-59}$ & $2.43\times10^{-17}$ & $4.72\times10^{-27}$ & 4.97 \\
5 & $6.38\times10^{-75}$ & $8.04\times10^{-49}$ & $7.01\times10^{-25}$ & 5.67 \\
6 & $9.17\times10^{-191}$ & $0.21\times10^{-160}$ & $3.47\times10^{-90}$ & 5.97 \\
\hline
\end{tabular}
\label{tab:early_errors}
\end{table}

The appendix experiment confirms that closer initial guesses accelerate the entrance into the asymptotic high-order regime, while even moderate or random starts still produce very rapid error decay within a few iterations. This observation is useful for documenting the numerical behavior of the underlying solver, but it should be interpreted separately from the main contribution of the present paper, which concerns prediction of $S_{\mathrm{mom}}$ across the parameter space.

\subsection{Comparison with Existing Parallel Schemes}
For completeness, we also report a compact comparison between the proposed scheme and two existing parallel methods from the literature. Table~\ref{tab:performance} summarizes the number of iterations to convergence, CPU time, memory usage, and maximum error, while Table~\ref{tab:order} reports theoretical and observed convergence orders.

\begin{table}[!htbp]
\centering
\caption{Computational performance of the proposed solver.}
\begin{tabular}{c c c c c}
\hline
Method & Iterations & CPU time (s) & Memory (MB) & Maximum error \\
\hline
Proposed Scheme~\cite{shams2026contractivity} & 5 & 0.012 & 123.65 & $2.05\times10^{-189}$\\
Existing Scheme~\cite{21} & 11 & 0.018 & 256.76 & $5.14\times10^{-45}$ \\
Existing Scheme~\cite{22} & 10 & 0.015 & 254.89 & $0.77\times10^{-37}$  \\
\hline
\end{tabular}
\label{tab:performance}
\end{table}

\begin{table}[!htbp]
\centering
\caption{Estimated convergence order of the methods.}
\begin{tabular}{c c c}
\hline
Method & Order (theoretical) & Order (observed) \\
\hline
Proposed Scheme~\cite{shams2026contractivity} & 6 & 5.97 \\
Existing Scheme~\cite{21} & 6 & 4.56 \\
Existing Scheme~\cite{22} & 6 & 3.19 \\
\hline
\end{tabular}
\label{tab:order}
\end{table}

These appendix results indicate that the proposed base scheme converges in fewer iterations and attains markedly smaller final errors than the two comparison methods, while preserving an observed order close to the theoretical sixth order. Again, this comparison is intended only as auxiliary solver validation. The main manuscript focuses instead on the AI-assisted diagnostic layer that predicts profile-based reliability over the $(\alpha,\beta)$ domain.

\section*{Funding}
The work reported in Sections~3 and~4 was supported by the Russian Science Foundation (grant no.~22-11-00055-P, \url{https://rscf.ru/en/project/22-11-00055/}, accessed on 10 June 2025) (Andrei Velichko). The material reported in Sections~1 and~5 was supported by the European Regional Development and Cohesion Funds (ERDF) 2021--2027 under Project AI4AM - EFRE1052 (Bruno Carpentieri). Bruno Carpentieri is also a member of the \textit{Gruppo Nazionale per il Calcolo Scientifico} (GNCS) of the \textit{Istituto Nazionale di Alta Matematica} (INdAM).

\section*{Author contributions}
Conceptualization: Bruno Carpentieri, Andrei Velichko, Mudassir Shams; Methodology: Bruno Carpentieri, Andrei Velichko, Mudassir Shams; Software: Andrei Velichko, Mudassir Shams; Validation: Bruno Carpentieri, Andrei Velichko, Mudassir Shams; Formal analysis: Bruno Carpentieri; Investigation: Bruno Carpentieri, Andrei Velichko, Mudassir Shams; Data curation: Andrei Velichko, Mudassir Shams; Writing---original draft: Bruno Carpentieri, Andrei Velichko, Mudassir Shams; Writing---review \& editing: Bruno Carpentieri, Andrei Velichko, Paola Lecca; Visualization: Andrei Velichko, Mudassir Shams; Supervision: Bruno Carpentieri, Andrei Velichko; Project administration: Bruno Carpentieri. All authors have read and agreed to the submitted version of the manuscript.


\section*{Disclosure of interest}
The authors report no conflict of interest.

\section*{Data availability statement}
All data products generated in this study (including grid-level metrics, raw and
smoothed proxy profiles, train/test split manifests, model predictions, and
aggregated summaries) are provided in the accompanying project directory. The
implementation scripts and detailed reproduction instructions will be made
publicly available in an online repository upon publication.


\nocite{cummings2025routing}
\bibliographystyle{unsrtnat}
\bibliography{references}

\end{document}